\documentclass[leqno,11pt]{amsart}
\usepackage{amsmath,amscd,amsthm,amsxtra}
\usepackage{epsfig,graphics}
\usepackage{amssymb,latexsym}
\usepackage{mathrsfs}
\usepackage{eucal}
\usepackage{hyperref}
\usepackage[all,cmtip]{xy}

\newtheorem{thm}{\bf Theorem}[section]
\newtheorem{df}[thm]{\bf Definition}
\newtheorem{prop}[thm]{\bf Proposition}
\newtheorem{cor}[thm]{\bf Corollary}
\newtheorem{lem}[thm]{\bf Lemma}
\newtheorem{rem}[thm]{\bf Remark}
\newtheorem{ex}[thm]{\bf Example}

\newcommand{\A}{\mathbb{A}}

\newcommand{\cP}{\mathscr{P}}

\newcommand{\pf}{\noindent{\bfseries Proof. }}
\newcommand{\ov}{\overline}

\newcommand{\M}{{\mathcal{M}}}

\newcommand{\gl}{\mathfrak{gl}}
\newcommand{\Z}{\mathbb{Z}}
\newcommand{\C}{\mathbb{C}}
\newcommand{\h}{\mathfrak{h}}
\newcommand{\te}{\widetilde{e}}
\newcommand{\tf}{\widetilde{f}}

\newcommand{\td}{\widetilde}
\newcommand{\mc}{\mathcal}
\newcommand{\mf}{\mathfrak}

\newcommand{\U}{{U}}
\newcommand{\K}{{K}}

\numberwithin{equation}{section}

\begin{document}
\title[ ]
{Crystal bases of $q$-deformed Kac modules over the quantum superalgebra $U_q(\gl(m|n))$}
\author{JAE-HOON KWON}
\address{Department of Mathematics \\ University of Seoul   \\  Seoul 130-743, Korea }
\email{jhkwon@uos.ac.kr }

\thanks{This work was  supported by Basic Science Research Program through the National Research Foundation of Korea (NRF) 
funded by the Ministry of  Education, Science and Technology (No. 2011-0006735).}

\begin{abstract} We introduce the notion of a crystal base of a finite dimensional $q$-deformed Kac module  over the quantum superalgebra $U_q(\gl(m|n))$, and  prove its existence and uniqueness. In particular,  we obtain the crystal base of a finite dimensional irreducible $U_q(\gl(m|n))$-module with typical highest weight.  We also show that the crystal base of a $q$-deformed Kac module is compatible with that of its irreducible quotient $V(\lambda)$  given by Benkart, Kang and Kashiwara when $V(\lambda)$ is an irreducible polynomial representation.
\end{abstract}

\maketitle


\section{Introduction}
\subsection{Background}
Let $U_q(\frak{g})$ be the  quantized enveloping algebra associated with a symmetrizable Kac-Moody algebra $\frak{g}$.
The  Kashiwara 's crystal base theory \cite{Kash1}  is  a study of  a certain nice basis of a $U_q(\frak{g})$-module $M$ at $q=0$, which still contains important combinatorial information on $M$, and  for the last couple of decades it has been one of the most important and successful tools in understanding the structures of modules over  both $\mathfrak{g}$ and $U_q(\frak{g})$.

Let $\gl(m|n)$ be a general linear Lie superalgebra over the complex numbers and  $U_q(\gl(m|n))$ its quantized enveloping algebra \cite{Ya}. In \cite{BKK}, Benkart, Kang and Kashiwara developed the  crystal base theory for a certain category of $U_q(\gl(m|n))$-modules called $\mathcal{O}_{int}$, which includes the irreducible factors appearing in a tensor power of the natural representation, often referred to as irreducible polynomial representations.  An irreducible $U_q(\gl(m|n))$-module in $\mathcal{O}_{int}$  can be viewed as a natural super-analogue of a finite dimensional irreducible $U_q(\gl(m+n))$-module since it is isomorphic to an irreducible polynomial  representation tensored  by a one dimensional representation. They proved that an irreducible module in $\mathcal{O}_{int}$   has a unique crystal base, and showed that its associated crystal can be realized in terms of $(m|n)$-hook semistandard tableaux \cite{BR} when it is polynomial. 

The crystal base of an irreducible module in $\mathcal{O}_{int}$ has several interesting features, which are not parallel to those over symmetrizable Kac-Moody (super) algebras \cite{J, Kash1}. For example, it is twisted compared to the usual crystal base over $U_q(\gl(m)\oplus\gl(n))$ in the sense that  it is a lower  crystal base as a $U_q(\gl(m|0))$-module  but  is a  upper  crystal base as a $U_q(\gl(0|n))$-module.
Its crystal also has fake highest weight elements so that it becomes quite non-trivial to prove that the crystal is connected.

The existence of a crystal base of an arbitrary finite dimensional irreducible  $U_q(\gl(m|n))$-module remains unknown unlike $U_q(\gl(m+n))$-modules, and one of the main obstacles in this problem is that a finite dimensional $U_q(\gl(m|n))$-module is not necessarily semisimple.  Indeed, the semisimplicity of a tensor power of the natural representation together with its polarizability is an essential ingredient in proving the existence of a crystal base of an irreducible polynomial representation. 

\subsection{Main results}
Let $P$ be the integral weight lattice of $\gl(m|n)$ and let $P^+$ be the set of integral weights  dominant with respect to its even subalgebra $\gl(m|n)_0=\gl(m|0)\oplus \gl(0|n)$. The finite dimensional  irreducible $\gl(m|n)$-modules with weights in $P$ are highest weight modules whose highest weights are parametrized by $P^+$. There is another important class of finite dimensional  $\gl(m|n)$-modules called {\it Kac modules}, which are indecomposable highest weight modules   parametrized by $P^+$ \cite{Kac2}. As the generalized or parabolic Verma modules do for $\gl(m+n)$-modules in a parabolic BGG category, they play the same role in a Kazhdan-Lusztig type  character formula for  finite dimensional irreducible $\gl(m|n)$-modules \cite{Br,CL,Se}. 

The purpose of this paper is to develop the crystal base theory of a $q$-deformed Kac module $K(\lambda)$ over $U_q(\gl(m|n))$ for $\lambda\in P^+$.
We introduce first the notion of a crystal base of $K(\lambda)$. Since $K(\lambda)$ does not necessarily belong to ${\mc O}_{int}$, we define modified Kashiwara operators $\te_0$ and $\tf_0$ on $K(\lambda)$ associated to the odd isotropic simple root $\alpha_0$, which are analogous to those on $U^-_q({\mf g})$ for a  symmetrizable Kac-Moody algebra ${\mf g}$ \cite{Kash1}. Then we show that $K(\lambda)$ has  a  unique crystal base  (Theorems \ref{main result} and \ref{main result - uniqueness}), which is the main result in this paper. In particular,  we obtain the crystal base of a finite dimensional irreducible  $U_q(\gl(m|n))$-module with typical highest weight. 

The key idea of the proof is to realize $K(\lambda)$ as a $q$-deformation of the exterior algebra $\Lambda((\mathbb{C}^m)^*\otimes \mathbb{C}^n)$ \cite{U} tensored by  an irreducible highest weight $U_q(\gl(m|n)_0)$-module $V_{m,n}(\lambda)$ with highest weight $\lambda$. Note that the $q$-deformation of the exterior algebra $\Lambda((\mathbb{C}^m)^*\otimes \mathbb{C}^n)$ has natural commuting actions of $U_q(\gl(m|0))$ and $ U_q(\gl(0|n))$ (or an action of $U_q(\gl(m|n)_0)$) and here we extend it to a $U_q(\gl(m|n))$-module by using a PBW type basis of $U_q(\gl(m|n))$. The crystal of $K(\lambda)$ has a simple description. As an underlying set, it is given by
\begin{equation}\label{crystal model}
\mathscr{P}(\Phi^-_1)\times \mathscr{B}_{m,n}(\lambda),
\end{equation}
where $\mathscr{P}(\Phi^-_1)$ is the power set of the set of negative odd roots of $\gl(m|n)$ and $\mathscr{B}_{m,n}(\lambda)$ is the crystal of $V_{m,n}(\lambda)$. Also, the crystal structure on \eqref{crystal model} can be described easily (Section \ref{combinatorics of Kac module}). 

We next show that the crystal base of $K(\lambda)$ is compatible with that of its irreducible quotient $V(\lambda)$ when $V(\lambda)$ is an irreducible polynomial representation (Theorem \ref{main result - compatibility}), that is, the canonical projection from $K(\lambda)$ to $V(\lambda)$ sends the crystal base of $K(\lambda)$ onto that of $V(\lambda)$. 
Hence we may regard the crystal of $V(\lambda)$ as a subgraph of the crystal of $K(\lambda)$. We give a combinatorial description of  its embedding (Section \ref{compatibility}) using the $(m|n)$-hook tableaux crystal model for $V(\lambda)$ and  the skew dual RSK algorithm introduced by Sagan and Stanley \cite{SS}. 
 

The paper is organized as follows. In Section \ref{preliminary}, we give necessary background on the quantum superalgebra $U_q(\gl(m|n))$. In Section \ref{crystal base theory}, we review the crystal base theory developed in \cite{BKK}. In Section \ref{crystal base of Kac modules}, we define the notion of a crystal base of a Kac module and state the main results, whose proofs are given in the following two sections. 
 
 \vskip 3mm
{\bf Acknowledgement} 
Part of this work was done while the author was visiting the Institute of Mathematics in Academia Sinica, Taiwan on January 2012. He would like to thank S.-J. Cheng for the invitation and helpful discussion, and the staffs for their warm hospitality.

\section{Quantum superalgebra $U_q(\gl(m|n))$}\label{preliminary}

\subsection{Lie superalgebra $\gl(m|n)$}
For non-negative integers $m,n$, let $[m|n]$ be a $\Z_2$-graded set with $[m|n]_0=\{\,\ov{m},\ldots,\ov{1}\,\}$, $[m|n]_1=\{\,1,\ldots,n\,\}$ and a linear ordering $\ov{m}<\ldots<\ov{1}<1<\ldots<n$. We denote by $|a|$ the degree of $a\in [m|n]$. 
Let $\C^{[m|n]}$ be the complex superspace
with a basis $\{\,v_a\,\vert\, a\in [m|n]\,\}$, where the parity of $v_a$ is $|a|$. 

Let
$\gl(m|n)$ denote the Lie
superalgebra of $[m|n]\times [m|n]$ complex matrices, which is spanned by  $E_{ab}$ $(a,b\in [m|n])$  with 1 in the $a$th row and the $b$th column, and $0$ elsewhere \cite{Kac1}. 

Let $P^{\vee} = \bigoplus_{a\in [m|n]}
\mathbb{Z} E_{aa}$ be the dual weight lattice  and  $\h=\C\otimes_{\Z} P^\vee$ the Cartan subalgebra of $\gl(m|n)$.
Define $\epsilon_{a}\in \h^*$ by $\langle E_{bb}, \epsilon_{a} \rangle
= \delta_{ab}$ for $a,b\in [m|n]$, where $\langle\cdot,\cdot\rangle$ denotes the natural pairing
on ${\mf h}\times {\mf h}^*$. Let $P=\bigoplus_{a\in [m|n]}
\mathbb{Z} \epsilon_{a}$ be the weight lattice of $\gl(m|n)$. For $\lambda=\sum_{a\in [m|n]}\lambda_a\epsilon_a\in P$, the parity of $\lambda$ is defined to be
$\lambda_1+\cdots+\lambda_n \mod 2$ and denoted by $|\lambda|$. Let $(\,\cdot\,|\,\cdot\,)$ denote a
symmetric bilinear form on $\frak{h}^* =
\mathbb{C} \otimes_{\mathbb{Z}} P$ given by $(\epsilon_a\,|\,\epsilon_b)=(-1)^{|a|}\delta_{ab}$  for $a,b\in [m|n]$.

Let $I=I_{m|n}=\{\,\ov{m-1},\ldots,\ov{1},0,1\ldots,n-1\,\}$, where we assume that $I_{m|0}=\{\,\ov{m-1},\ldots,\ov{1}\,\}$ and $I_{0|n}=\{\,1\ldots,n-1\,\}$. Then with respect to the Borel subalgebra spanned by $E_{ab}$ ($a\le b$), the set of simple roots of $\gl(m|n)$ is $\Pi=\{\,\alpha_k\,|\,k\in I_{m|n}\,\}$, where
\begin{equation*}
\alpha_k=
\begin{cases}
\epsilon_{\ov{i+1}}-\epsilon_{\ov{i}}, & \text{if $k=\ov{i}\in I_{m|0}$},\\
\epsilon_{\ov{1}}-\epsilon_{1}, & \text{if $k=0$},\\
\epsilon_j-\epsilon_{j+1}, & \text{if $k=j\in I_{0|n}$}.
\end{cases}
\end{equation*}
Note that $(\alpha_k|\alpha_k)=2$ (resp. $-2$) for $k\in I_{m|0}$ (resp. $I_{0|n}$), and $(\alpha_0|\alpha_0)=0$.
Let $Q=\bigoplus_{k\in I}\Z\alpha_k$ be the root lattice, and $Q^\pm=\pm\sum_{k\in I}\Z_{\geq 0}\alpha_k$.  We have a partial ordering on $P$, where $\lambda\geq \mu$ if and only if $\lambda-\mu\in Q^+$ for $\lambda, \mu\in P$.

The set of positive roots, even positive roots and odd positive roots are given by
\begin{equation*}
\begin{split}
\Phi^+&=\{\, \epsilon_a-\epsilon_b \,|\,   a<b  \,\}, \\
\Phi^+_0&=\{\,\epsilon_a-\epsilon_b\,|\,a<b,  \ |a|= |b|\,\}=\{\,\alpha\in \Phi^+\,|\,(\alpha|\alpha)=\pm2\,\},\\
\Phi^+_1&=\{\,\epsilon_a-\epsilon_b\,|\,a<b,  \ |a|\neq |b|\,\}=\{\,\alpha\in \Phi^+\,|\,(\alpha|\alpha)=0\,\},\\
\end{split}
\end{equation*}
respectively. 

The simple coroot $h_k$ ($k\in I$)
corresponding to $\alpha_k$ is the unique element in $P^{\vee}$
satisfying $l_k \langle h_k, \lambda \rangle = (\alpha_k | \lambda)$  for all $\lambda \in P$, where $l_k =1$ (resp. $l_k=-1$) for $k\in I_{m|0}\cup\{ 0\}$ 
(resp. $k\in I_{0|n}$). 

Let 
\begin{equation*}
P^+=\left\{\,\lambda=\sum_{a\in [m|n]}\lambda_a\epsilon_a \in P\ \Bigg|\ \lambda_{\ov{m}}\geq\ldots\geq\lambda_{\ov{1}},\ \ \lambda_1\geq \ldots\geq \lambda_n\,\right\},
\end{equation*}
which is the set of dominant integral weights for $\gl(m|0)\oplus\gl(0|n)\subset\gl(m|n)$.
For $\lambda\in P^+$, let 
\begin{equation*}
\begin{split}
\lambda_{+}=\sum_{i=1}^m\lambda_{\ov{i}}\epsilon_{\ov{i}},\ \ \ \ \
\lambda_{-}=\sum_{j=1}^n\lambda_{j}\epsilon_{j},
\end{split}
\end{equation*}
Also, we let $\delta=\epsilon_{\ov{m}}+\cdots+\epsilon_{\ov{1}}-\epsilon_1-\cdots-\epsilon_n$,
where we have $\delta_+=\epsilon_{\ov{m}}+\cdots+\epsilon_{\ov{1}}$ and $\delta_-=-\epsilon_1-\cdots-\epsilon_n$.

\subsection{Quantum superalgebra $U_q(\gl(m|n))$} We assume that $q$ is an indeterminate.
The  quanum superalgebra $U_q(\gl(m|n))$ is the
associative superalgebra (or $\Z_2$-graded algebra) over $\mathbb{Q}(q)$ generated by
$e_k$, $f_k$ $(\,k\in I\,)$ and $q^h$ $(\,h\in P^{\vee}\,)$, which are
subject to the following relations \cite{BKK,Ya}:
{\allowdisplaybreaks
\begin{equation*}
\begin{split}
& {\rm deg}(q^h)=0,\ \  {\rm deg}(e_k)={\rm deg}(f_k)=|\alpha_k|,\\
& q^0=1, \quad q^{h +h'}=q^{h}q^{h'}, \quad \ \ \ \ \ \ \ \ \ \ \ \ \ \ \ \ \ \ \ \ \ \ \ \  \ \ \ \ \ \ \ \ \text{for} \ \
h, h' \in P^{\vee}, \\ &q^h e_k=q^{\langle h,\alpha_k\rangle}
e_k q^h, \quad q^h f_k=q^{-\langle h,\alpha_k\rangle} f_k q^h, \\ 
& e_k
f_l-(-1)^{|\alpha_k||\alpha_l|}f_l e_k =\delta_{kl}\frac{t_k-t_k^{-1}}
{q_k-q^{-1}_k}, \\ & e_k e_l - (-1)^{|\alpha_k||\alpha_l|} e_l e_k = f_k f_l
- (-1)^{|\alpha_k||\alpha_l|} f_l f_k =0, \quad \!\text{if $(\alpha_k | \alpha_l) =0$},
\\ & e_k^2 e_l-(q+q^{-1}) e_k e_l e_k+e_l e_k^2= 0, \quad\ \ \ \ \ \ \ \ \  \ \ \ \ \ \ \ \ \ \ \ \ \ 
\text{if $k\neq 0$ and $(\alpha_k | \alpha_l) \neq 0$},
\\ & f_k^2 f_l-(q+q^{-1}) f_k f_l f_k+f_l f_k^2= 0, \quad\ \ \ \ \ \ \ \ \ \ \ \ \ \ \ \ \  \ \ \ \  \text{if $k\neq 0$ and $(\alpha_k | \alpha_l) \neq 0$},
\\ & e_0 e_{\ov{1}} e_0 e_{1} + e_{\ov{1}} e_0
e_{1} e_{0} + e_{0} e_{1} e_{0} e_{\ov{1}}  + e_{1} e_{0} e_{\ov{1}} e_{0} -(q+q^{-1}) e_{0} e_{\ov{1}}
e_{1} e_{0} =0, \\ & f_0 f_{\ov{1}} f_0 f_{1} + f_{\ov{1}} f_0
f_{1} f_{0} + f_{0} f_{1} f_{0} f_{\ov{1}}  + f_{\bar 1} f_{0} f_{\ov{1}} f_{0} -(q+q^{-1}) f_{0} f_{\ov{1}}
f_{1} f_{0} =0.
\end{split}
\end{equation*}
}Here, $q_k = q^{l_k}$ and  $t_k = q^{l_k h_k}$.

For simplicity, we will assume the following notations throughout the paper:
\begin{itemize}
\item[$\cdot$]  $\U=U_q(\gl(m|n))$,

\item[$\cdot$]  $\U^\pm$ : the subalgebras  generated by $e_k$ and $f_k$ ($k\in I$), respectively,

\item[$\cdot$] $\U^0$ : the subalgebra generated by $q^h$ ($h\in P^\vee$),

\item[$\cdot$] $\U_{m,n}$ : the subalgebra  generated by $q^h, e_k, f_k$ ($h\in P^\vee, k\in I\setminus\{0\}$),

\item[$\cdot$] $\U^\pm_{m,n} = \U_{m,n}\cap \U^\pm$,

\item[$\cdot$]  $\U_{m|0}$ : the subalgebra generated by $q^{E_{\ov{i}\,\ov{i}}}, e_k, f_k$ ($i=1,\ldots,m, k\in I_{m|0}$),

\item[$\cdot$]  $\U_{0|n}$ : the subalgebra generated by $q^{E_{jj}}, e_k, f_k$ ($j=1,\ldots,n, k\in I_{0|n}$).
\end{itemize}

There is a Hopf superalgebra structure on $U$, where the comultiplication $\Delta$ is given by
\begin{equation*}
\begin{split}
\Delta(q^h)&=q^h\otimes q^h, \\ \Delta(e_k)&=e_k\otimes t_k^{-1} + 1\otimes e_k, \\  
\Delta(f_k)&=f_k\otimes 1+ t_k\otimes f_k,
\end{split}
\end{equation*}
the antipode $S$ is given by 
\begin{equation*}
S(q^h)=q^{-h}, \ \ S(e_k)=-e_kt_k, \ \  S(f_k)=-t_k^{-1}f_k,
\end{equation*}
and the couint $\varepsilon$ is given by  $\varepsilon(q^h)=1$, $\varepsilon(e_k)=\varepsilon(f_k)=0$ for $h\in P^\vee$ and  $k\in I$.

A $U$-module $M$ is always understood to be a $U$-supermodule, that is, $M=M_0\oplus M_1$ with $U_iM_j\subset M_{i+j}$ for $i,j\in \Z_2$ (see \cite{Kac1} for basic notions related with superalgebras). 
We also have a superalgebra structure on $U\otimes U$ with the multiplication $(u_1\otimes u_2)(v_1\otimes v_2)=(-1)^{|u_2||v_1|}(u_1v_1)\otimes (u_2v_2)$, where $|u|$ denotes the degree of a homogeneous element $u\in U$. Hence, we have a $U$-module structure on $M_1\otimes M_2$ via the comultiplication $\Delta$ for $U$-modules $M_1$ and $M_2$. 

For $\mu\in P$, $M_\mu=\{\,m\,|\,q^h m =q^{\langle h,\mu \rangle}m\ (h\in P^\vee)\,\}$ is called a weight space of $M$ with weight $\mu$. When $M=\bigoplus_{\mu\in P}M_\mu$, we say that $M$ has a weight space decomposition. Throughout this paper, we assume that the $\Z_2$-grading on $M$ is induced from the parity of its weights when $M$ has a weight space decomposition.
 
Note that $U_{m,n}\cong  U_{m|0}\otimes U_{0|n}$ as a $\mathbb{Q}(q)$-algebra, and $\U_{m|0}$ (resp. $U_{0|n}$) is isomorphic to the quantized enveloping algebra  $U_q(\gl_m)$ (resp. $U_{q^{-1}}(\gl_n)$), whose $\mathfrak{sl}_2$-copy is  generated by $e_k$, $f_k$ and $t_k$ for $k\in I_{m|0}$ (resp. $k\in I_{0|n}$). But, when we consider a $U_{0|n}$-module in this paper, we understand $\U_{0|n}$ as $U_v(\gl_n)$, whose $\mathfrak{sl}_2$-copy is  generated by ${\mf e}_k=e_k$, ${\mf f}_k=f_k$ and ${\mf t}_k=t_k^{-1}$ with $v=v_k=q_k^{-1}=q$ for $k\in I_{0|n}$. We denote by $(\,\cdot\,|\,\cdot\,)'=-(\,\cdot\,|\,\cdot\,)$ the symmetric bilinear form on the weight lattice of $U_{0|n}$ so that $(\alpha_k|\alpha_k)'=2$ for $ k\in I_{0|n}$.

For $\lambda\in P^+$, let $V_{m,n}(\lambda)$ be an irreducible $U_{m,n}$-module with highest weight $\lambda$, and let $V_{m|0}(\lambda_+)$ (resp. $V_{0|n}(\lambda_-)$) the irreducible highest weight module over $U_{m|0}$ (resp. $U_{0|n}$) with highest weight $\lambda_+$ (resp. $\lambda_-$). Note that $V_{m,n}(\lambda)\cong V_{m|0}(\lambda_+)\otimes V_{0|n}(\lambda_-)$ as a $U_{m,n}$-module. 

\subsection{PBW type basis}
We have $\U^\pm =\bigoplus_{\alpha\in Q^\pm}\U^{\pm}_\alpha$, where $\U^{\pm}_\alpha=\{\,u\,|\,q^h u q^{-h}=q^{\langle h,\alpha\rangle}u\ (h\in P^\vee)\,\}$. For $x\in \U^+_\alpha$, $y\in \U^+_\beta$, we define the super $q$-bracket by
\begin{equation*}
[x,y]_q=xy-(-1)^{|\alpha||\beta|}q^{-(\alpha|\beta)}yx.
\end{equation*}
For $\alpha\in \Phi^+$ with $\alpha=\alpha_{k}+\alpha_{k+1}+\ldots+\alpha_l$ ($k<l$), we define
\begin{equation*}
e_\alpha=[[\cdots[[e_k,e_{k+1}]_q,e_{k+2}]_q\cdots]_q,e_{l}]_q.
\end{equation*}
Here, we assume a linear ordering $\ov{m-1}<\ldots<\ov{1}<0<1<\ldots<n-1$ on $I$, and $k+1$ denotes the adjacent element in $I$, which is larger than $k\in I$.

We define a linear ordering on $\Phi^+$ by
\begin{equation*}
\alpha < \beta  \ \ \, \Longleftrightarrow \ \ \text{$(a<c)$ or ($a=c$ and $b<d$)},
\end{equation*}
for $\alpha,\beta \in \Phi^+$ with $\alpha=\epsilon_a-\epsilon_b$ and $\beta=\epsilon_{c}-\epsilon_d$.
For $\alpha, \beta\in \Phi^+$ with $\alpha< \beta$, it is straightforward to check the following commutation relations (cf.\cite{Ya});
\begin{equation}\label{commutation relation}
\begin{split}
[e_{\alpha},e_\beta]_q=
\begin{cases}
(-1)^{|\alpha||\delta|}(q^{(\alpha|\gamma-\alpha)}-q^{-(\alpha|\gamma-\alpha)})e_{\delta}e_{\gamma}, & \text{for $a<c<b<d$},\\
e_{\gamma}, & \text{for $b=c$},\\
0, & \text{otherwise},
\end{cases}
\end{split}
\end{equation}
where we assume that $\alpha=\epsilon_a-\epsilon_b$, $\beta=\epsilon_c-\epsilon_d$, $\gamma=\epsilon_a-\epsilon_d$ and $\delta=\epsilon_c-\epsilon_b$. In particular, we have $e_\alpha^2=0$  for $\alpha\in\Phi^+_1$.

\begin{prop}[Proposition 10.4.1 in \cite{Ya}]\label{Yamane}
Let
\begin{equation*}
{B}^+=\Big\{\,\overrightarrow{\prod_{\alpha\in\Phi^+}}e_\alpha^{m_\alpha} \,\Big|\, \text{$m_{\alpha}\in\Z_{\geq 0}$ for $|\alpha|=0$ and $m_{\alpha}=0,1$ for $|\alpha|=1$}\,\Big\},
\end{equation*} 
where the product is taken in the order of $<$ on $\Phi^+$. Then
${B}^+$ is a $\mathbb{Q}(q)$-basis of $\U^+$.
\end{prop}

Let $\sharp$ be the $\mathbb{Q}(q)$-linear anti-involution on $\U$ given by $e_k^\sharp=f_k$, $f_k^\sharp=e_k$ and $(q^h)^\sharp=q^{h}$ for $k\in I$ and $h\in P^\vee$. Then $B^-=({B}^+)^\sharp$ is a $\mathbb{Q}(q)$-basis of $\U^-$, and  
\begin{equation}\label{triangular decomposition}
\U\cong \U^-\otimes \U^0\otimes \U^+
\end{equation}
as a  $\mathbb{Q}(q)$-vector space with a basis $\{\,u^- q^h u^+\,|\,u^{\mp}\in {B}^{\mp},\ h\in P^\vee\,\}$ \cite[Theorem 10.5.1]{Ya}.  
Since $\U$ is a Hopf superalgebra, we have a $\mathbb{Q}(q)$-algebra homomorphism  ${\rm ad} : \U \longrightarrow {\rm End}_{\mathbb{C}(q)}(\U)$ given by  
\begin{equation}\label{adjoint action}
\begin{split}
{\rm ad}(q^h)(u)&=q^h u q^{-h}, \\
{\rm ad}(e_k)(u)&=(e_k u -(-1)^{|\alpha_k||u|}u e_k) t_k, \\
{\rm ad}(f_k)(u)&=f_k u -(-1)^{|\alpha_k||u|}t_k u t_k^{-1} f_k, \\
\end{split}
\end{equation}
for $h\in P^\vee$, $k\in I$ and a homogeneous element $u$.
For $\alpha\in \Phi^+$, we put 
$f_\alpha = e_{\alpha}^\sharp$. If $\alpha=\alpha_{k}+\alpha_{k+1}+\ldots+\alpha_l$ ($k<l$), then we have
\begin{equation*}
f_\alpha = {\rm ad}(f_{l})\circ \cdots \circ {\rm ad}(f_{k+2})\circ {\rm ad}(f_{k+1})(f_k).
\end{equation*}
By applying $\sharp$ to \eqref{commutation relation}, we have
\begin{equation}\label{commutation relation for f}
\begin{split}
[f_{\beta},f_\alpha]_q=
\begin{cases}
(-1)^{|\alpha||\delta|}(q^{(\alpha|\gamma-\alpha)}-q^{-(\alpha|\gamma-\alpha)})f_{\gamma}f_{\delta}, & \text{for $a<c<b<d$},\\
f_{\gamma}, & \text{for $b=c$},\\
0, & \text{otherwise}.
\end{cases}
\end{split}
\end{equation}

\subsection{Parabolic decomposition and a $q$-deformed wedge space}
For $\alpha\in \Phi^+_1$ with $\alpha=\epsilon_{\ov{i}}-\epsilon_j$, we put
\begin{equation}\label{f_S}
\begin{split}
{\bf f}_\alpha &= {\rm ad}(f_{j-1})\circ \cdots \circ {\rm ad}(f_{1}) \circ {\rm ad}(f_{\ov{i-1}})\circ \cdots \circ  {\rm ad}(f_{\ov{1}})(f_0) \\
&= {\rm ad}(f_{\ov{i-1}})\circ \cdots \circ {\rm ad}(f_{\ov{1}}) \circ {\rm ad}(f_{j-1})\circ \cdots \circ {\rm ad}(f_{1})(f_0).
\end{split}
\end{equation}
Let $\K$ be the subalgebra of $\U^-$ generated by ${\bf f}_\alpha$ $(\alpha\in \Phi_1^+)$.

Let us define a linear ordering on $\Phi^+_1$ by
\begin{equation*}
\alpha\prec \beta  \ \ \, \Longleftrightarrow \ \ \text{$(b<d)$ or ($b=d$ and $a>c$)}
\end{equation*}
for $\alpha,\beta \in \Phi^+_1$ with $\alpha=\epsilon_a-\epsilon_b$ and $\beta=\epsilon_{c}-\epsilon_d$. For $S\subset \Phi^+_1$ with $S=\{\,\beta_1\prec \cdots \prec \beta_r\,\}$, we put
\begin{equation}\label{monomial F_S}
{\bf f}_S = {\bf f}_{\beta_1}\cdots {\bf f}_{\beta_r}.
\end{equation}  
Here we assume that ${\bf f}_S=1$ when $S=\emptyset$. 
It is straightforward to check that for $\alpha=\epsilon_{\ov{i}}-\epsilon_j$, $\beta=\epsilon_{\ov{k}}-\epsilon_l \in \Phi^+_1$ with $\alpha\prec \beta$
\begin{equation}\label{commutation relation 2}
\begin{split}
{\bf f}_\alpha  {\bf f}_\beta &=- q{\bf f}_\beta{\bf f}_\alpha \ \ \ \ \ \ \ \ \ \ \ \ \ \ \ \ \ \ \ \ \ \ \ \ \ \ \ \, \text{for \ ($i=k$,\ $j<l$) or ($i<k$,\  $j=l$)},\\
{\bf f}_\alpha  {\bf f}_\beta &=- {\bf f}_\beta{\bf f}_\alpha,\ \ \ \ \ \ \ \ \ \ \ \ \ \ \ \ \ \ \ \ \ \ \ \ \ \ \ \, \text{for \ $i>k$,\ $j<l$}, \\
{\bf f}_\alpha  {\bf f}_\beta &=- {\bf f}_\beta{\bf f}_\alpha +(q-q^{-1}){\bf f}_\gamma{\bf f}_\delta,  \  \ \ \ \ \ \  \text{for \ $i<k$,\ $j<l$}, \\
{\bf f}_\alpha^2&=0,
\end{split}
\end{equation}
where $\gamma=\epsilon_{\ov{i}}-\epsilon_l$ and $\delta=\epsilon_{\ov{k}}-\epsilon_j$. 

\begin{lem}\label{KL decomposition}
${B}_{\K}=\{\,{\bf f}_S\,|\,S\subset \Phi^+_1\,\}$ is a $\mathbb{Q}(q)$-basis of $K$, and 
\begin{equation*}
\U^- \cong \K\otimes \U^-_{m,n}
\end{equation*}
as a $\mathbb{Q}(q)$-vector space.
\end{lem}
\pf Given $u\in U^-$, we have $u=u_1u_2$ for some $u_1\in K$ and $u_2\in U^-_{m,n}$ by \eqref{commutation relation for f} and \eqref{f_S}. Hence  $B=\{\,u_1u_2 \,|\, u_1\in {B}_{\K}, \ u_2\in {B}^-\cap\U_{m,n}^-\}$  spans $U^-$ since $B_K$ spans $K$ by \eqref{commutation relation 2} and $B^-\cap U^-_{m,n}$ is a basis of $U^-_{m,n}$ by Proposition \ref{Yamane}. Considering the dimension of $\U^-_\alpha$ for each $\alpha\in Q^-$, we see that $B$ is linearly independent, and hence a $\mathbb{Q}(q)$-basis of $\U^-$. In particular, $B_K$ is a $\mathbb{Q}(q)$-basis of $K$.  \qed
\vskip 2mm

Let us define another linear orderings on $\Phi^+_1$ by
\begin{equation*}
\begin{split}
&\alpha \prec' \beta \ \ \, \Longleftrightarrow \ \ \text{$(a>c)$ or $(a=c, \ b>d)$}, \\
&\alpha \prec'' \beta \ \, \, \Longleftrightarrow \ \ \text{$(a>c)$ or $(a=c, \ b<d)$}, \\
\end{split}
\end{equation*}
for $\alpha, \beta\in \Phi^+_1$ with $\alpha=\epsilon_a-\epsilon_b$ and $\beta=\epsilon_c-\epsilon_d$.
Then as in \eqref{monomial F_S}, we may define ${\bf f}'_S$, ${\bf f}''_S$ ($S\subset \Phi^+_1$) and ${B}'_{\K}$, ${B}''_{\K}$ under $\prec'$, $\prec''$, respectively. Suppose that ${\bf f}_S\in U^-_{\alpha}$ with
\begin{equation}\label{weight of S}
\alpha=-\sum_{i=1}^ma_i\epsilon_{\ov{i}}+\sum_{j=1}^n b_j\epsilon_j
\end{equation}
for some $a_i$ and $b_j\in \Z_{\geq 0}$. By the first two relations in \eqref{commutation relation 2}, we have 
\begin{equation}\label{relation f and f'}
{\bf f}_S  = \pm {\bf f}''_S= \pm(-q)^{\sum_{i}a_i(a_i-1)/2}{\bf f}'_S.
\end{equation}
In particular, $B'_{\K}$ and $B''_{\K}$ are $\mathbb{Q}(q)$-bases of $\K$. 
Note that
\begin{equation}\label{norm of alpha}
(\alpha|\alpha)=\sum_{i}a_i^2-\sum_{j}b_j^2=\sum_{i}a_i(a_i+1)-\sum_{j}b_j(b_j+1).
\end{equation}

\begin{rem}{\rm 
Note that $\K$ is a $q$-deformed exterior algebra generated by $(\mathbb{Q}(q)^m)^*\otimes \mathbb{Q}(q)^n$. This can be explained as follows. For $\alpha=\epsilon_{\ov{i}}-\epsilon_j\in\Phi^+_1$, let $k(\alpha)=i+m(j-1)\in\mathbb{N}$. We have a bijection from the set of non-empty $S\subset \Phi^+_1$ to the set of non-empty subsets in $\{\,1,\ldots,mn\,\}$, by sending $S=\{\,\beta_1, \ldots ,\beta_r\,\}$ to $\{\,k(\beta_1), \ldots ,k(\beta_r)\,\}$. 
Consider the $q$-deformed wedge space $\Lambda=\bigoplus_{r\geq 1}\Lambda^r$ in \cite[Section 3.3]{U} (here we replace $n$ and $l$ in \cite{U} by $m$ and $n$, respectively). 
Define a map 
\begin{equation*}
\kappa : \K \longrightarrow \Lambda
\end{equation*}
by $\kappa({\bf f}_S)=(-1)^{\sum b_j(b_j+1)/2}u_{k(\beta_1)}\wedge \cdots \wedge u_{k(\beta_r)}$ for $S=\{\,\beta_1, \ldots ,\beta_r\,\}\subset \Phi^+_1$ with $b_j$  as in \eqref{weight of S}. By comparing \eqref{commutation relation 2} and \cite[Proposition 3.16]{U} (see also \cite[Remark 3.17 (ii)]{U}), we see that $\kappa$ is an injective $\mathbb{Q}(q)$-algebra homomorphism.
}
\end{rem}

\subsection{Kac modules}
Let $L$ be the subalgebra of $\U$ generated by $\U_{m,n}$ and $e_0$. For $\lambda\in P^+$, we extend $V_{m,n}(\lambda)$ to an $L$-module in an obvious way, and  define 
\begin{equation*}
K(\lambda) = \U\otimes_{L}V_{m,n}(\lambda).
\end{equation*}
to be the induced $\U$-module. Since $L\cong \U^-_{m,n}\otimes \U^0 \otimes \U^+$ as a $\mathbb{Q}(q)$-vector space,  
\begin{equation}\label{K linear iso}
K(\lambda) \cong \K\otimes V_{m,n}(\lambda),
\end{equation}
as a $\mathbb{Q}(q)$-vector space by \eqref{triangular decomposition} and Lemma \ref{KL decomposition}. Note that $K(\lambda)=\bigoplus_{\mu\leq \lambda}K(\lambda)_\mu$ with ${\rm dim}K(\lambda)_\lambda=1$, and $K(\lambda)=U^- 1_\lambda$, where  $K(\lambda)_\lambda=\mathbb{Q}(q)1_\lambda$.
Hence $K(\lambda)$ is a finite dimensional highest weight module with highest weight $\lambda$. We call $K(\lambda)$ a ($q$-deformed) {\it Kac module with highest weight $\lambda$}.  We define $V(\lambda)$ to be the  maximal irreducible quotient of $K(\lambda)$, and denote the image of $1_\lambda$ by $v_\lambda$. 

\subsection{Classical limits and typicality of Kac modules}
Let us consider classical limit of $U$-modules. We leave the detailed verification to the reader since the argument below are almost identical to the case of symmetrizable Kac-Moody algebras (see \cite{Ja,Lu}). 

Let ${\bf A}=\mathbb{Q}[q,q^{-1}]$. Let $M$ be a highest weight $U$-module generated by a highest weight vector $u$ of weight $\lambda\in P$. 
We define
\begin{equation*}
M_{\bf A}=\sum_{r\geq 0,\, i_1,\ldots,i_r\in I}{\bf A}f_{i_1}\ldots f_{i_r}u,\ \ M_{\mu,{\bf A}}=\sum_{\substack{r\geq 0,\,i_1,\ldots,i_r\in I \\ \lambda-\alpha_{i_1}-\cdots-\alpha_{i_r}=\mu}}{\bf A}f_{i_1}\ldots f_{i_r}u.
\end{equation*}
Then $M_{{\bf A}}=\bigoplus_{\mu}M_{\mu,{\bf A}}$, and 
${\rm rank}_{\bf A}M_{\mu,{\bf A}}=\dim_{\mathbb{Q}(q)}M_\mu$
since ${\bf A}$ is a principal ideal domain and $M_{\mu,{\bf A}}$ is a torsion free (hence free) ${\bf A}$-module.
It is easy to check that the ${\bf A}$-module $M_{\bf A}$ is invariant under $e_k$, $f_k$, $q^{h}$ and $\frac{q^h-q^{-h}}{q-q^{-1}}$ for $k\in I$ and $h\in P^\vee$.

Let $\varphi : {\bf A} \longrightarrow \mathbb{C}$ be a $\mathbb{Q}$-algebra homomorphism given by $\phi(f(x))=f(1)$. Set
\begin{equation*}
\ov{M}=M_{\bf A}\otimes_{\bf A}\mathbb{C}, \ \  \ov{M}_{\mu}=M_{\mu, {\bf A}}\otimes_{\bf A}\mathbb{C}.
\end{equation*}
Here $\mathbb{C}$ is understood to be an ${\bf A}$-module via $\varphi$. We have $\ov{M}=\bigoplus_{\mu}\ov{M}_\mu$ with $\dim_{\mathbb{C}}\ov{M}_\mu={\rm rank}_{\bf A}M_{\mu,{\bf A}}$.
Let $\ov{e_k}$, $\ov{f_k}$ and $\ov{h}$ be the $\mathbb{C}$-linear endomorphisms on $\ov{M}$ induced from $e_k$, $f_k$ and $\frac{q^h-q^{-h}}{q-q^{-1}}$ for $k\in I$ and $h\in P^\vee$. Let $\ov{U}_{\ov{M}}$ be the subalgebra of ${\rm End}_{\mathbb{C}}(\ov{M})$ generated by $\ov{e_k}$, $\ov{f_k}$ and $\ov{h}$ for $k\in I$ and $h\in P^\vee$. Then there exists a $\mathbb{C}$-algebra homomorphism from $U(\gl(m|n))$ to $\ov{U}_{\ov{M}}$ sending $E_k$,  $F_{k}$ and $h$ to $\ov{e_k}$, $\ov{f_k}$ and $\ov{h}$, where $U(\gl(m|n))$ is the enveloping algebra of $\gl(m|n)$, and $E_k$, $F_k$ denote the root vectors in $\gl(m|n)$ corresponding to $\alpha_k$, $-\alpha_k$, respectively.
Hence, $\ov{M}$ is a $U(\gl(m|n))$-module. 

We see that $\ov{K(\lambda)}$ is a Kac module over $\gl(m|n)$ when $M=K(\lambda)$  for $\lambda\in P^+$ by comparing the dimensions of weight spaces of $K(\lambda)$ and $\ov{K(\lambda)}$.
For $\lambda\in P^+$, $\lambda$ is called {\it typical} if $\langle \alpha, \lambda+\rho\rangle\neq 0$ for all $\alpha\in \Phi^+_1$, where $\rho$ is the Weyl vector of $\gl(m|n)$ given by $\rho=\frac{1}{2}\sum_{\alpha\in \Phi^+_0}\alpha-\frac{1}{2}\sum_{\beta\in \Phi^+_1}\beta$. It is shown by Kac \cite[Proposition 2.9]{Kac2} that $\ov{K(\lambda)}$ is irreducible when $\lambda$ is typical. 

\begin{prop}
For typical  $\lambda\in P^+$, $K(\lambda)$ is irreducible.
\end{prop}
\pf Suppose that $K(\lambda)$ is not irreducible. Let $M$ be a proper submodule of $K(\lambda)$. We may assume that $M$ is a highest weight module. From the above argument, it follows that $\ov{M}$ is a proper submodule of $\ov{K(\lambda)}$, which is a contradiction since $\lambda$ is typical. Hence $K(\lambda)$ is irreducible. \qed

\section{Review on crystal base theory}\label{crystal base theory}
\subsection{Crystal bases for $U_q(\gl_{m})$-modules}\label{review on crystal}
Let us briefly recall the notion of crystal bases for integrable $U_q(\gl_m)$-modules (we refer the reader to \cite{Kash1} for more details in a general setting). Here we assume that $U_q(\gl_m)=U_{m|0}$. 

Let $M$ be an integrable $U_q(\gl_m)$-module with weight decomposition $M=\bigoplus_{\lambda\in P}M_\lambda$ and ${\rm dim}M_\lambda<\infty$. For $u\in M_{\lambda}$ and $k\in I_{m|0}$, we have
\begin{equation*}
u=\sum_{r\geq 0, -\langle h_k,\lambda \rangle}f_k^{(r)}u_r,
\end{equation*}
where $e_ku_r=0$ for all $r\ge 0$. 
Here,
\begin{equation*}
[r]_k =\frac{q_k^r-q_k^{-r}}{q_k-q_k^{-1}},\quad
[r]_k!=\prod_{s=1}^{r}[s]_k, \quad
f_k^{(r)}=\frac{1}{[r]_k!}f_k^r.
\end{equation*}
Then the Kashiwara operators  are defined by
\begin{equation}\label{Kashiwara operator}
\tilde{e}_k u=
\sum_{r \ge 1} f_k^{(r-1)}u_r, \ \ \   
\tilde{f}_k u=
\sum_{r \ge 0} f_k^{(r+1)}u_r.
\end{equation}

Let $\mathbb{A}$ denote the subring of $\mathbb{Q}(q)$ consisting
of all rational functions which are regular at $q=0$. A pair $(L,B)$ is called a {\it lower} {\it crystal base of $M$} if  
\begin{itemize}
\item[(1)] $L$ is an $\A$-lattice of $M$, where $L=\bigoplus_{\lambda\in P}L_{\lambda}$ with $L_{\lambda}=L\cap M_{\lambda}$,
\item[(2)]  $\tilde{e}_k L\subset L$ and $\tilde{f}_k L\subset L$ for $k\in I_{m|0}$,
\item[(3)] $B$ is a $\mathbb{Q}$-basis of $L/qL$, where $B=\bigsqcup_{\lambda\in P}B_{\lambda}$ with
$B_{\lambda}=B\cap(L/qL)_{\lambda}$,

\item[(4)] $\tilde{e}_kB \subset B\sqcup \{0\}$,
$\tilde{f}_k B\subset B\sqcup \{0\}$ for $k\in I_{m|0}$,

\item[(5)] for $b,b'\in B$ and $k\in I_{m|0}$,  
$\tilde{f}_k b = b'$ if and only if $b=\tilde{e}_k b'$.
\end{itemize}
For $b\in B$ and $k\in I_{m|0}$, we set
\begin{equation}\label{epsilon phi}
\varepsilon_k(b)= \max\{\,r\in\mathbb{Z}_{\geq
0}\,|\,\tilde{e}_k^r b\neq 0\,\}, \ \ \ \varphi_k(b)=
\max\{\,r\in\mathbb{Z}_{\geq 0}\,|\,\tilde{f}_k^r b\neq 0\,\}.
\end{equation}
Let $k=\mathbb{Q}(q^{\frac{1}{2}})$ and let $\psi_M : k\otimes_{\mathbb{Q}(q)}M \longrightarrow k\otimes_{\mathbb{Q}(q)}M$ be a $\mathbb{Q}(q)$-linear isomorphism given by $\psi_M(u)=q^{-\frac{(\lambda|\lambda)}{2}}u$ for $u\in M_\lambda$. For $k\in I_{m|0}$, we define $\te_k^{up}, \tf_k^{up} : M \longrightarrow M$ by 
\begin{equation*}
\te_k^{up}=\psi_M\circ \te_k\circ\psi_M^{-1}, \ \ \ \tf_k^{up}=\psi_M\circ \tf_k\circ\psi_M^{-1}.
\end{equation*}
Note that $\te_k^{up}$ and $\tf_k^{up}$ are well-defined $\mathbb{Q}(q)$-linear operators on $M$ even when $(\lambda|\lambda)\not\in 2\Z$ for some weight $\lambda$ of $M$. Then a pair $(L,B)$ is  called a {\it upper crystal base of $M$} if it satisfies the above conditions (1)--(5) with respect to $\te_k^{up}$ and $\tf_k^{up}$. 

\begin{rem}{\rm  If $\psi_M(M)= M$ and $(L,B)$ is a lower crystal base of $M$, then   $( \psi_M(L),  \psi_M(B))$ is a upper crystal base of $M$. Also, if $M$ is an integrable highest weight module with highest weight $\lambda$ with a lower crystal base $(L,B)$, then $q^{\frac{(\lambda|\lambda)}{2}}\psi_M(M)=M$ and $(q^{\frac{(\lambda|\lambda)}{2}}\psi_M(L),q^{\frac{(\lambda|\lambda)}{2}}\psi_M(B))$ is a upper crystal base of $M$. }
\end{rem}

Let $M_i$  ($i=1,2$) be  integrable $U_q(\gl_m)$-modules with lower crystal bases $(L_i,B_i)$.
Then $(L_1\otimes L_2, B_1\otimes B_2)$ is a lower crystal base of $M_1\otimes M_2$. The operators $\te_k$, $\tf_k$ act on $B_1\otimes B_2$;
\begin{equation}\label{lower tensor product rule}
\begin{split}
&\tilde{e}_k(b_1\otimes b_2)= \begin{cases}
\tilde{e}_k b_1 \otimes b_2, & \text{if $\varphi_k(b_1)\geq\varepsilon_k(b_2)$}, \\ 
b_1 \otimes \tilde{e}_k b_2, & \text{if $\varphi_k(b_1)<\varepsilon_k(b_2)$},\\
\end{cases}
\\
&\tilde{f}_k(b_1\otimes b_2)=
\begin{cases}
\tilde{f}_k b_1 \otimes b_2, & \text{if $\varphi_k(b_1)>\varepsilon_k(b_2)$}, \\
 b_1 \otimes \tilde{f}_k b_2, & \text{if $\varphi_k(b_1)\leq\varepsilon_k(b_2)$}. 
\end{cases}
\end{split}
\end{equation}
On the other hand, let $(L_i^{up}, B_i^{up})$ be upper crystal bases of $M_i$ ($i=1,2$), respectively.
Denote by $M_1\otimes_+ M_2 $ be a tensor product with a $U_q(\gl_m)$-module structure induced from the comultiplication 
\begin{equation}\label{Delta_+}
\begin{split}
\Delta_+(q^h)&=q^h\otimes q^h, \\ \Delta_+(e_k)&=e_k\otimes 1 + t_k\otimes e_k, \\  
\Delta_+(f_k)&=f_k\otimes t_k^{-1}+ 1\otimes f_k,
\end{split}
\end{equation}
for $h\in P^\vee$ and  $k\in I_{m|0}$.
Then $(L^{up}_1\otimes L^{up}_2, B^{up}_1\otimes B^{up}_2)$ is a upper crystal base of $M_1\otimes_+ M_2$. The operators $\te^{up}_k$, $\tf^{up}_k$ act on $B_1\otimes B_2$ in the same way as in \eqref{lower tensor product rule}.
 
\subsection{Crystal bases of $U$-modules in $\mathcal{O}_{int}$}\label{crystal base for Oint}
Let $\mathcal{O}_{int}$ be  the category of $U$-modules $M$ satisfying the following conditions:
\begin{itemize}
\item[(1)] $M=\bigoplus_{\lambda\in P}M_\lambda$ with ${\rm dim}M_\lambda<\infty$,

\item[(2)] $M$ is an integrable $U_{m,n}$-module,

\item[(3)]  if $M_{\lambda} \neq 0$, then
$\langle h_0,\lambda \rangle \geq 0$,

\item[(4)] if $\langle h_0,\lambda \rangle= 0$, then
$e_0 M_{\lambda}=f_0 M_{\lambda}=0$.
\end{itemize}

Let us review the notion of crystal bases for $U$-modules in $\mathcal{O}_{int}$ \cite{BKK}.
Let $M=\bigoplus_{\lambda \in P} M_{\lambda}\in \mathcal{O}_{int}$ and let $u\in M_\lambda$ be given.
For $k\in I_{m|0}$, we define $\te_ku$ and $\tf_ku$ to be as in \eqref{Kashiwara operator} with $u$ as an element of a $U_{m|0}$-module. For $k\in I_{0|n}$, we define $\te_ku$ and $\tf_ku$ to be $\tilde{\mf e}^{up}_k u$ and $\tilde{\mf f}^{up}_k u$ with $u$ as an element of an $U_{0|n}$-module (see Section \ref{review on crystal}). 
For $k=0$, we define 
\begin{equation}\label{e0 f0 for Oint}
\tilde{e}_{0} u = q^{-1} t_0 e_0 u, \quad \text{and} \quad
\tilde{f}_{0} u = f_0 u.
\end{equation}
Then a pair $(L,B)$ is called a {\it crystal base of $M$} if 
\begin{itemize}
\item[(1)] $L$ is an $\A$-lattice of $M$, where  $L=\bigoplus_{\lambda\in P}L_{\lambda}$ with $L_{\lambda}=L\cap
M_{\lambda}$,
\item[(2)] $\tilde{e}_k L\subset L$ and $\tilde{f}_k L\subset L$
for $k\in I$,

\item[(3)] $B$ is a pseudo-basis of $L/qL$ (i.e. 
$B=B^{\bullet}\cup(-B^{\bullet})$ for a $\mathbb{Q}$-basis
$B^{\bullet}$ of $L/qL$),

\item[(4)] $B=\bigsqcup_{\lambda\in P}B_{\lambda}$ with
$B_{\lambda}=B\cap(L/qL)_{\lambda}$,

\item[(5)] $\tilde{e}_kB \subset B\sqcup \{0\}$,
$\tilde{f}_k B\subset B\sqcup \{0\}$ for $k\in I$,

\item[(6)] for $b,b'\in B$ and $k\in I$,
$\tilde{f}_k b = b'$ if and only if $b=\tilde{e}_k b'$.
\end{itemize}
The set $B / \{\pm 1 \}$ has an $I$-colored oriented graph
structure, where  $b\stackrel{k}{\rightarrow}
b'$ if and only if $\tilde{f}_k b=b'$ for $k\in I$ and  $\,b, b' \in B / \{\pm
1\}$. We call $B / \{\pm 1 \}$ the {\it crystal} of $M$.

For $b\in B$ and $k\in I$, let $\varepsilon_k(b)$ and $\varphi_k(b)$ be as in \eqref{epsilon phi}.
Let us recall the {tensor product
rule} for the crystal bases of $U$-modules in
$\mathcal{O}_{int}$ (see \cite[Proposition 2.8]{BKK}). Let $M_i$ $(\,i=1,2\,)$ be a
$U$-module in $\mathcal{O}_{int}$
with a crystal base $(L_i,B_i)$. Then  $(L_1\otimes L_2,B_1\otimes B_2)$ is a crystal
base of $M_1\otimes M_2$.  For $k\in I_{m|0}$,  $\te_k$ and $\tf_k$ act on $B_1\otimes B_2$ as in \eqref{lower tensor product rule}. Note that $\Delta({\mf e}_k)={\mf e}_k\otimes {\mf t}_k + 1\otimes {\mf e}_k$ and $\Delta({\mf f}_k)={\mf f}_k\otimes 1 + {\mf t}_k^{-1}\otimes {\mf f}_k$ for $k\in I_{0|n}$, and they coincide with $\Delta_+({\mf e}_k)$ and $\Delta_+({\mf f}_k)$ (see \eqref{Delta_+}) if we exchange the order of the tensor product. Hence, for $k\in I_{0|n}$, the formulas for $\te_k$ and $\tf_k$ on $B_1\otimes B_2$ are given by exchanging the positions of  tensor factors in \eqref{lower tensor product rule} since $(L_i,B_i)$  are upper crystal bases of $M_i$ as $U_{0|n}$-modules, that is,{\allowdisplaybreaks
\begin{equation}\label{upper tensor product rule}
\begin{split}
&\tilde{e}_k(b_1\otimes b_2)= \begin{cases}
 b_1 \otimes \tilde{e}_kb_2, & \text{if $\varphi_k(b_2)\geq\varepsilon_k(b_1)$}, \\ 
\tilde{e}_k b_1 \otimes b_2, & \text{if $\varphi_k(b_2)<\varepsilon_k(b_1)$},\\
\end{cases}
\\
&\tilde{f}_k(b_1\otimes b_2)=
\begin{cases}
b_1 \otimes \tilde{f}_k b_2, & \text{if $\varphi_k(b_2)>\varepsilon_k(b_1)$}, \\
\tilde{f}_k b_1 \otimes  b_2, & \text{if $\varphi_k(b_2)\leq\varepsilon_k(b_1)$}. 
\end{cases}
\end{split}
\end{equation}}
For $k=0$, we have{\allowdisplaybreaks
\begin{equation}\label{tensor product rule for 0}
\begin{split}
\tilde{e}_0(b_1\otimes b_2)=&
\begin{cases}
\tilde{e}_0 b_1\otimes b_2, & \text{if }\langle h_0,{\rm
wt}(b_1)\rangle>0, \\ \pm b_1\otimes \tilde{e}_0 b_2, & \text{if
}\langle h_0,{\rm wt}(b_1)\rangle=0,
\end{cases}
\\
\tilde{f}_0(b_1\otimes b_2)=&
\begin{cases}
\tilde{f}_0 b_1\otimes b_2, & \text{if }\langle h_0,{\rm
wt}(b_1)\rangle>0, \\ \pm b_1\otimes \tilde{f}_0 b_2, & \text{if
}\langle h_0,{\rm wt}(b_1)\rangle=0.
\end{cases}
\end{split}
\end{equation}}
Here, ${\rm wt}$ denotes the weight function and the $\pm$ sign  depends on the parity of ${\rm wt}(b_1)$.

\subsection{Semistandard tableaux}
Let us recall some basic background on tableaux (see \cite{BR,Ful}), which will be used in later sections.  

Let $\mc{P}$ be the set of partitions.
A partition $\lambda=(\lambda_i)_{i\geq 1}$ is identified with a Young diagram. We denote by $\lambda'=(\lambda_i')_{i\geq 1}$ its conjugate.

Let $\A$ be a linearly
ordered  set with a $\mathbb{Z}_2$-grading $\A=\A_0\sqcup\A_1$. For a skew  Young diagram  $\lambda/\mu$, a tableau $T$ obtained by
filling $\lambda/\mu$ with entries in $\A$ is called
$\A$-semistandard  if (1) the entries in each row (resp. column) are
weakly increasing from left to right (resp. from top to bottom), (2)
the entries in $\A_0$ (resp. $\A_1$) are strictly increasing in each
column (resp. row).  We say that $\lambda/\mu$ is the shape of
$T$, and write ${\rm sh}(T)=\lambda/\mu$. 
We denote by ${ SST}_{\A}(\lambda/\mu)$ the set of all
$\A$-semistandard tableaux of shape $\lambda/\mu$.

Let $\lambda\in\mathcal{P}$ be given. For $T\in SST_\A(\lambda)$ and
$a\in \A$, we denote by $a \rightarrow T$ the tableau obtained from $T$ by applying the usual Schensted column insertion (see \cite{Ful} and \cite{BR} for its super analogue).
We also need the following variation of the Schensted's column insertion. Let  $\lambda^\pi$ be the skew Young
diagram obtained by $180^{\circ}$-rotation of $\lambda$ (sometimes called of anti-normal shape).
For $T\in SST_\A(\lambda^\pi)$ and
$a\in \A$, we define $T \leftarrow a$ to be the tableau of an
anti-normal shape obtained from $T$ by applying the following
procedure; 
\begin{itemize}
\item[(1)] If $|a|=0$, then let $a'$ be the largest entry in the right-most
column which is smaller than or equal to $a$. If $|a|=1$, then  let $a'$ be the largest entry, which is smaller than $a$.

\item[(2)]  replace {$a'$} by
{$a$}. If there is no such $a'$, put {$a$} at the top of the column
and stop the procedure, 

\item[(3)] repeat (1) and (2) on the next column
with {$a'$}.
\end{itemize}
For a finite word $w=w_1\ldots w_r$ with letters in $\A$, we define $(w \rightarrow T)=(
w_n\rightarrow(\cdots(w_1\rightarrow T)))$ and $(T\leftarrow w)=(\cdots((T\leftarrow{w_r})\leftarrow{w_{r-1}}\,)
\cdots)\leftarrow{w_1}$.

\subsection{Crystal bases of polynomial representations}\label{crystal of poly}
Let us review the results on the crystal bases of irreducible polynomial representations of $U$ \cite{BKK}. 
Let $V=\bigoplus_{b\in [m|n]}\mathbb{Q}(q)v_b$ be the natural representation of $U$. Then $V$ has a crystal base $(\mathscr{L},\mathscr{B})$ where $\mathscr{L}=\bigoplus_{b\in [m|n]}\A v_b$ and $\mathscr{B}=\{\,\pm\ov{v_b}\,|\,b\in [m|n]\,\}$.  For simplicity, we identify $\mathscr{B}/\{\pm1\}$ with $[m|n]$ as a set, and the crystal of $V$ is given by
\begin{equation*}
\ov{m}\ \stackrel{^{\ov{m-1}}}{\longrightarrow}\ \ov{m-1}\ \stackrel{^{\ov{m-2}}}{\longrightarrow}
\cdots\stackrel{^{\ov{1}}}{\longrightarrow}\ \ov{1}\ \stackrel{^0}{\longrightarrow}\ 1 \
\stackrel{^1}{\longrightarrow}\cdots\stackrel{^{n-2}}{\longrightarrow} n-1\stackrel{^{n-1}}{\longrightarrow} n.
\end{equation*}

For $r\geq 1$, $(\mathscr{L}^{\otimes r},\mathscr{B}^{\otimes r})$ is a crystal base of $V^{\otimes r}$. Let $\mathscr{W}$ be the set of all finite words with the
letters in $[m|n]$. The empty word is denoted by $\emptyset$.
Then $\mathscr{W}$ is a crystal of the tensor algebra  
since we may identify each non-empty word $w=w_1\cdots w_r$ with
$w_1\otimes\cdots\otimes w_r\in \mathscr{B}^{\otimes r}/\{\pm1\}$, where
$\{\emptyset\}$ forms a trivial crystal of weight $0$.

Let 
\begin{equation*}
\td{P}^+=\left\{\,\lambda=\sum_{a\in [m|n]}\lambda_a\epsilon_a \in P\ \Bigg|\ \lambda_{\ov{m}}\geq\ldots\geq\lambda_{\ov{1}}\geq \lambda'_1\geq \lambda'_2\geq \ldots\,\right\}.
\end{equation*}
For $r\geq 1$, $V^{\otimes r}$ is completely reducible and each irreducible $U$-module in $V^{\otimes r}$, which we call an irreducible polynomial representation, is isomorphic to $V(\lambda)$ for some $\lambda\in \td{P}^+$ with $\sum_{a\in [m|n]}\lambda_a=r$.

Let $\mc{P}_{m|n}=\{\,\mu=(\mu_i)_{i\geq 1}\in\mc{P}\,|\,\mu_{m+1}\leq n\,\}$ which is called the set of $(m|n)$-hook partitions. 
Then the map sending $\mu=(\mu_i)_{i\geq 1}$ to  $\mu_1\epsilon_{\ov{m}}+\cdots +\mu_1\epsilon_{\ov{1}}+\nu'_1\epsilon_1+\cdots+\nu'_n\epsilon_n$ is a bijection from $\mc{P}_{m|n}$ to $\td{P}^+$, where $\nu=(\nu_i)_{i\geq 1}$ is given by $\nu_i=\mu_{m+i}$.

Now, let $\lambda^\circ$ be the Young diagram (or partition) corresponding to $\lambda\in \td{P}^+$.
For $T\in SST_{\mathscr{B}}(\lambda^\circ)$, let $T(i,j)$ denote  the entry of $T$ located in
the $i$-th row from the top and the $j$-th column from the left.
Then we choose an embedding 
\begin{equation}\label{admissible reading}
\psi : SST_{\mathscr{B}}(\lambda^\circ) \rightarrow \mathscr{W}
\end{equation}
 by reading the entries of $T$ in $SST_{\mathscr{B}}(\lambda^\circ)$
in such a way that $T(i,j)$ should be read before $T(i+1,j)$ and
$T(i,j-1)$ for each $i,j$. The image of
$SST_{\mathscr{B}}(\lambda^\circ)$ under $\psi$ together with $0$ is stable under
$\te_k,\tf_k$ ($k\in I$) and  the induced  $I$-colored oriented graph structure does not depend on the choice of $\psi$ \cite[Theorem 4.4]{BKK}. Moreover, $SST_{\mathscr{B}}(\lambda^\circ)$ is connected with a unique highest weight element \cite[Theorem 4.8]{BKK}.

\begin{thm}\label{BKK}{\rm(\cite[Theorem 5.1]{BKK})} 
For $\lambda\in \td{P}^+$, $V(\lambda)$ has a unique crystal base  $(\mathscr{L}(\lambda),\mathscr{B}(\lambda))$ such that $\mathscr{L}(\lambda)_\lambda=\A v_\lambda$, and  $\mathscr{B}(\lambda)/\{\pm1\} \cong SST_{\mathscr{B}}(\lambda^\circ)$, that is, there is a weight preserving isomorphism of $I$-colored oriented graphs.  
\end{thm}

\section{Crystal bases of Kac modules}\label{crystal base of Kac modules}
\subsection{Crystal base of $K(\lambda)$ as a $U_{m,n}$-module} 
Let us consider the $U_{m,n}$-action on $K(\lambda)$  ($\lambda\in P^+$).
For simplicity, let us write $u\cdot v = {\rm ad}(u)(v)$ for $u,v\in \U$ (see  \eqref{adjoint action}).  

Let $k\in I\setminus\{0\}$ be given. 
By \eqref{f_S}, we have for $\alpha\in\Phi^+_1$
\begin{equation}\label{sl2 action on basis of K}
e_k\cdot{\bf f}_\alpha= {\bf f}_{\alpha-\alpha_k}, \ \ \   
f_k\cdot{\bf f}_\alpha= {\bf f}_{\alpha+\alpha_k}, \ \ \ t_k\cdot {\bf f}_\alpha = q^{-(\alpha|\alpha_k)}{\bf f}_\alpha.
\end{equation}
Here we assume that ${\bf f}_{\alpha\pm\alpha_k}=0$ when $\alpha\pm\alpha_k\not\in\Phi^+_1$.
For  $S=\{\,\beta_1\prec \cdots \prec \beta_r\,\}\subset \Phi^+_1$ and $v\in V_{m,n}(\lambda)$,  we have {\allowdisplaybreaks
\begin{equation}\label{action on K} 
\begin{split}
e_k({\bf f}_S\otimes v)&=(e_k{\bf f}_{\beta_1}\cdots {\bf f}_{\beta_r})\otimes_{L} v \\
&=\sum_{i=1}^r{\bf f}_{\beta_1}\cdots {\bf f}_{\beta_{i-1}}(e_k\cdot{\bf f}_{\beta_i})(t_k^{-1}\cdot{\bf f}_{\beta_{i+1}})\cdots (t_k^{-1}\cdot{\bf f}_{\beta_{r}})\otimes_{L} t_k^{-1}v\\
&\ \ \ \   + {\bf f}_{\beta_1}\cdots {\bf f}_{\beta_r}\otimes_{L} e_k v,\\
f_k({\bf f}_S\otimes v)&=(f_k{\bf f}_{\beta_1}\cdots {\bf f}_{\beta_r})\otimes_{L} v \\
&=\sum_{i=1}^r(t_k\cdot{\bf f}_{\beta_1})\cdots (t_k\cdot {\bf f}_{\beta_{i-1}})(f_k\cdot{\bf f}_{\beta_i}){\bf f}_{\beta_{i+1}}\cdots {\bf f}_{\beta_{r}}\otimes_{L} v_0 \\
&\ \ \ \  + (t_k\cdot{\bf f}_{\beta_1})\cdots (t_k\cdot {\bf f}_{\beta_{r}})\otimes_{L} f_k v,\\
t_k({\bf f}_S\otimes v)&=(t_k{\bf f}_{\beta_1}\cdots {\bf f}_{\beta_r})\otimes_{L} v 
= (t_k\cdot{\bf f}_{\beta_1})\cdots (t_k\cdot {\bf f}_{\beta_{r}})\otimes_{L} t_k v.
\end{split}
\end{equation}}

\begin{prop} For $\lambda\in P^+$, we have  as a $\U_{m,n}$-module
$$K(\lambda) \cong K(0)\otimes V_{m,n}(\lambda).$$
\end{prop}
\pf By Lemma \ref{KL decomposition}, $\{\,{\bf f}_S\otimes u_0\,|\,{\bf f}_S\in B_K\,\}$ is a $\mathbb{Q}(q)$-basis of $K(0)$, where $u_0$ is the highest weight vector in $V_{m,n}(0)$. Define a map $\phi : K(\lambda) \longrightarrow K(0)\otimes_{\mathbb{Q}(q)} V_{m,n}(\lambda)$ by $\phi({\bf f}_S\otimes_{L} v)=({\bf f}_S\otimes_{L} u_0)\otimes_{\mathbb{Q}(q)} v$ for  ${\bf f}_S\in {B}_{\K}$ and $v\in V_{m,n}(\lambda)$, which is a well defined $\mathbb{Q}(q)$-linear isomorphism by \eqref{K linear iso}. Then $\phi$ is a $\U_{m,n}$-module homomorphism by \eqref{action on K}, and hence an isomorphism.
\qed\vskip 2mm

Set
\begin{equation*}
\begin{split}
\mathscr{L}(K)&=\bigoplus_{S\subset \Phi^+_1}\A\,{\bf f}_S {1}_0\subset K(0 ),\\
\mathscr{B}(K)&=\left\{\  {\bf f}_S {1}_0 \!\! \mod q\mathscr{L}(K) \ \big|\ S\subset \Phi^+_1\,\right\}.
\end{split}
\end{equation*}

\begin{prop}\label{lower crystal base of K(0)}  $(\mathscr{L}(K),\mathscr{B}(K))$ is a lower crystal base of $K(0)$ as a $U_{m|0}$-module.
\end{prop}
\pf  For $j=1,\ldots,n$, put 
\begin{equation*}
\begin{split}
\Phi^+_1[j]&=\{\,\epsilon_{\ov{i}}-\epsilon_j \,|\,i=1,\ldots,m\,\},\\
K(0){[j]}&=\bigoplus_{S\subset \Phi^+_1[j]}\mathbb{Q}(q){\bf f}_S 1_0.
\end{split}
\end{equation*}
By \eqref{sl2 action on basis of K} and \eqref{action on K}, it is straightforward to check that as a $U_{m|0}$-module
\begin{equation*}
K(0){[j]}\cong V_{m|0}(0)\oplus\bigoplus_{i=1}^mV_{m|0}(-\epsilon_{\ov{1}}-\cdots-\epsilon_{\ov{i}}).
\end{equation*}
Since ${\bf f}_S={\bf f}_{S{[1]}}\cdots {\bf f}_{S{[n]}}$ with $S{[j]}=S\cap \Phi^+_1[j]$ for $S\subset \Phi^+_1$, the map 
\begin{equation}\label{phi m|0}
\phi_{m|0} : K(0)\longrightarrow K(0){[1]}\otimes\cdots\otimes K(0){[n]}
\end{equation}
given by $\phi_{m|0}({\bf f}_S 1_0)={\bf f}_{S{[1]}}1_0\otimes \cdots \otimes{\bf f}_{S{[n]}}1_0$  is an isomorphism of $U_{m|0}$-modules by \eqref{action on K}. Therefore, $(\mathscr{L}(K),\mathscr{B}(K))$ is a lower crystal base of $K(0)$ as a $U_{m|0}$-module since
$(\mathscr{L}(K){[j]},\mathscr{B}(K){[j]})$ is a lower crystal base of $K(0){[j]}$ and $\phi_{m|0}(\mathscr{L}(K))=\mathscr{L}(K){[1]}\otimes \cdots\otimes \mathscr{L}(K){[n]}$, where
\begin{equation*}
\begin{split}
\mathscr{L}(K){[j]}&=\bigoplus_{S\subset \Phi^+_1[j]}\A\,{\bf f}_S {1}_0,\\
\mathscr{B}(K){[j]}&=\left\{\ {\bf f}_S {1}_0\!\!\! \mod q\mathscr{L}(K){[j]} \ \big|\ S\subset \Phi^+_1[j]\,\right\}.
\end{split}
\end{equation*}
\qed\vskip 3mm

Next, we set
\begin{equation*}
\begin{split}
\mathscr{L}(K)'&=\bigoplus_{S\subset \Phi^+_1}\A\,q^{\omega(S)}{\bf f}'_S {1}_0\subset K(0 ),\\
\mathscr{B}(K)'&=\left\{\ q^{\omega(S)}{\bf f}'_S {1}_0\!\!\! \mod q\mathscr{L}(K)'\ \Big|\ S\subset \Phi^+_1\,\right\}.
\end{split}
\end{equation*}
where  $\omega(S)=\sum_{i=1}^ma_i(a_i-1)/2$ for $S\subset \Phi^+_1$ with $\sum_{\beta\in S}\beta=-\sum_{i=1}^ma_i\epsilon_{\ov{i}}+\sum_{j=1}^n b_j\epsilon_j$.

\begin{prop}\label{upper crystal base of K(0)'}  $(\mathscr{L}(K)',\mathscr{B}(K)')$ is a upper crystal base of $K(0)$ as a $U_{0|n}$-module.
\end{prop}
\pf For $i=1,\ldots,m$, put 
\begin{equation*}
\begin{split}
\Phi^+_1\left[\ov{i}\right]&=\{\,\epsilon_{\ov{i}}-\epsilon_j \,|\, j=1,\ldots,n\,\},\\
K(0){\left[\ov{i}\right]}&=\bigoplus_{S\subset \Phi^+_1\left[\ov{i}\right]}\mathbb{Q}(q){\bf f}'_S 1_0 \subset K(0).
\end{split}
\end{equation*}
By \eqref{sl2 action on basis of K} and \eqref{action on K},
\begin{equation*}
K(0){\left[\ov{i}\right]}\cong V_{0|n}(0)\oplus \bigoplus_{j=1}^mV_{0|n}(\epsilon_{1}+\cdots+\epsilon_{j}),
\end{equation*}
 as a $U_{0|n}$-module. Since $K(0){\left[\ov{i}\right]}$ has a lower crystal base
\begin{equation*}
\begin{split}
\left(\bigoplus_{S\subset \Phi^+_1\left[\ov{i}\right]}\A\,{\bf f}'_S {1}_0,\ \
\left\{\  \ov{{\bf f}'_S {1}_0} \ \big|\ S\subset \Phi^+_1\left[\ov{i}\right]\,\right\} \right),
\end{split}
\end{equation*}
it has a upper crystal base
\begin{equation*}
\begin{split}
\mathscr{L}(K)'{\left[\ov{i}\right]}&=\bigoplus_{S\subset \Phi^+_1\left[\ov{i}\right]}\A\,q^{-\frac{(\beta_S|\beta_S)'}{2}}{\bf f}'_S {1}_0\\
\mathscr{B}(K)'{\left[\ov{i}\right]}&=\left\{\  q^{-\frac{(\beta_S|\beta_S)'}{2}}{\bf f}'_S {1}_0\!\!\! \mod q\mathscr{L}(K)'{\left[\ov{i}\right]}\ \big|\ S\subset \Phi^+_1\left[\ov{i}\right]\,\right\},
\end{split}
\end{equation*}
where $\beta_S=\sum_{\beta\in S}\beta$. Note that for $S\subset \Phi^+_1\left[\ov{i}\right]$, we have $-(\beta_S|\beta_S)'/2=(\beta_S|\beta_S)/2=a_i(a_i-1)/2$, where $a_i=(\epsilon_{\ov{i}}|\beta_S)$.

The map 
\begin{equation}\label{phi 0|n}
\phi_{0|n} : K(0)\longrightarrow K(0){[\ov{1}]}\otimes\cdots\otimes K(0){[\ov{m}]}
\end{equation}
given by $\phi_{m|0}({\bf f}'_S 1_0)={\bf f}'_{S{[\ov{1}]}}1_0\otimes \cdots \otimes{\bf f}'_{S{[\ov{m}]}}1_0$  is an isomorphism of $U_{0|n}$-modules  by \eqref{action on K} and the fact that ${\bf f}'_S={\bf f}'_{S{[\ov{1}]}}\cdots {\bf f}'_{S{[\ov{n}]}}$ with $S{\left[\ov{i}\right]}=S\cap \Phi^+_1\left[\ov{i}\right]$ for $S\subset \Phi^+_1$. Moreover, we have $\phi_{0|n}(\mathscr{L}(K)')=\mathscr{L}(K)'[\ov{1}]\otimes_+\cdots\otimes_+ \mathscr{L}(K)'[\ov{m}]$ since $\omega(S)=\omega(S[\ov{1}])+\cdots+\omega(S[\ov{m}])$. Therefore, $(\mathscr{L}(K)',\mathscr{B}(K)')$ is a upper crystal base of $K(0)$ (see (1.4.7) in \cite{Kash1}). Here we used $\otimes_+$ to emphasize that the comultiplication is $\Delta_+$ \eqref{Delta_+} with respect to ${\mf e}_k$, ${\mf f}_k$ and ${\mf t}_k$ ($k\in I_{0|n}$), where the order of tensor factors are reversed. \qed

\begin{prop}\label{lower = upper} 
We have $\mathscr{L}(K)=\mathscr{L}(K)'$ and $\mathscr{B}(K)/\{\pm1\}=\mathscr{B}(K)'/\{\pm1\}$.
\end{prop}
\pf It follows directly from \eqref{relation f and f'}.\qed\vskip 3mm


\subsection{Crystal base of $K(\lambda)$}
Let us define the notion of a crystal base of a Kac module $K(\lambda)$. The definition of  a crystal base in \cite{BKK} is not available for $K(\lambda)$ since it does not belong to $\mathcal{O}_{int}$ in general. So we give a different definition, which is based on \cite[Section 3]{Kash1}.

Let  $e'_0$ be  a $\mathbb{Q}(q)$-linear operator on $U^-$ characterized by
\begin{itemize}
\item[(1)] $e'_0(f_k) = \delta_{0 k}$ for $k\in I$,

\item[(2)] $e'_0(uv) = e'_0(u)v + (-1)^{|\alpha|}q^{(\alpha_0|\alpha)}u e'_0(v)$ for $u\in U^-_{\alpha}, v\in U^-$.
\end{itemize}
It is straightforward  to check that $e_0'$ is well-defined on $U^-$. We have another $\mathbb{Q}(q)$-linear operator on $U^-$ given by
\begin{itemize}
\item[(1)] $e''_0(f_k) = \delta_{0 k}$ for $k\in I$,

\item[(2)] $e''_0(uv) = e''_0(u)v + (-1)^{|\alpha|}q^{-(\alpha_0|\alpha)}u e''_0(v)$ for $u\in U^-_{\alpha}, v\in U^-$.
\end{itemize}
These two operators satisfy the following
\begin{equation*}
e_0P-(-1)^{|\alpha|}Pe_0= \frac{t_0 e''_0(P)-t_0^{-1}e'_0(P)}{q-q^{-1}}
\end{equation*}
for $P\in U^-_\alpha$ (see \cite[Section 3.3]{Kash1}).

\begin{lem} We have
\begin{equation*}
U^-={\rm Ker}\,e'_0 \oplus {\rm Im}\,f_0.
\end{equation*}
Here we understand $f_0$ as a linear operator acting on $U^-$ by the left multiplication. 
\end{lem}
\pf It is easy to see that $e'_0({\bf f}_\alpha) =0$ for $\alpha\in \Phi^+_0$. Suppose that $\alpha\in \Phi^+_1$.
If $\alpha=\alpha_0+\alpha_k$, then ${\bf f}_\alpha={\rm ad}(f_k)(f_{0})=f_k{f}_{0} - q^{-(\alpha_k|\alpha_0)}{f}_{0} f_k$ and $e'_0({\bf f}_\alpha)= q^{-(\alpha_k|\alpha_0)}f_k- q^{-(\alpha_k|\alpha_0)} f_k=0$. If  $\alpha=\beta+\alpha_k$ for some $\beta\in \Phi^+_1$ and $k\neq 0$, then ${\bf f}_\alpha={\rm ad}(f_k)({\bf f}_\beta)=f_k{\bf f}_\beta - q^{-(\alpha_k|\beta)}{\bf f}_\beta f_k$ and we have $e'_0({\bf f}_\alpha)=0$  by induction on the height of $\alpha$. Hence, we have for $\alpha\in \Phi^+$
\begin{equation}\label{e0 action on PBW}
e'_0({\bf f}_\alpha)=
\begin{cases}
1, & \text{if $\alpha=\alpha_0$}, \\
0, & \text{otherwise}.
\end{cases}
\end{equation}

Let $W_1$ be the $\mathbb{Q}(q)$-span of $B_1=\{\,{\bf f}_S u_0 \,|\,S\subset \Phi^+_1\, (\alpha_0\not\in S),\, u_0\in U_{m,n}\,\}$ and let $W_2$ be the $\mathbb{Q}(q)$-span of $B_2=\{\,{\bf f}_S u_0 \,|\,S\subset \Phi^+_1\, (\alpha_0\in S),\, u_0\in U_{m,n}\,\}$. Then $U^-=W_1\oplus W_2$. Since  $W_1\subset {\rm Ker}\,e'_0$, $W_2\subset {\rm Im}\,f_0$ by \eqref{e0 action on PBW}, and ${\rm Ker}\,e'_0 \cap {\rm Im}\,f_0=\{0\}$, we have $W_1= {\rm Ker}\,e'_0$ and  $W_2= {\rm Im}\,f_0$. \qed\vskip 3mm

For $\lambda\in P^+$, we may identify $K(\lambda)$ with $U^-/ I_\lambda$, where $I_\lambda$ is a left $U_-$-ideal generated by  $f_k^{\langle h_k, \lambda\rangle+1}$ ($k\in I\setminus \{0\}$). Since $e'_0(I_\lambda)=0$, it induces a $\mathbb{Q}(q)$-linear map on $K(\lambda)$, which we still denote by $e'_0$. 
For $u\in K(\lambda)$, we define
\begin{equation}\label{e0 f0 for K}
\te_0 u = e'_0(u), \ \ \ \tf_0 u = f_0 u.
\end{equation}

\begin{df}\label{def of crystal base of Kac module}{\rm 
For $\lambda\in P^+$, a pair $(L,B)$ is a {\it crystal base of $K(\lambda)$} if 
\begin{itemize}
\item[(1)] $L$ is an $\A$-lattice of $M$, where  $L=\bigoplus_{\mu\in P}L_{\mu}$ with $L_{\mu}=L\cap
K(\lambda)_{\mu}$,
\item[(2)] $\tilde{e}_k L\subset L$ and $\tilde{f}_k L\subset L$
for $k\in I$,

\item[(3)] $B$ is a pseudo-basis of $L/qL$, where  $B=\bigsqcup_{\mu\in P}B_{\mu}$ with
$B_{\mu}=B\cap(L/qL)_{\mu}$,

\item[(4)] $\tilde{e}_kB \subset B\sqcup \{0\}$,
$\tilde{f}_k B\subset B\sqcup \{0\}$ for $k\in I$,

\item[(5)] for any $b,b'\in B$ and $k\in I$, we have
$\tilde{f}_k b = b'$ if and only if $b=\tilde{e}_k b'$.
\end{itemize}
Let us call the $I$-colored oriented graph $B/\{\pm1\}$ a {\it crystal of $K(\lambda)$}.
}
\end{df}

\subsection{Main results} Now, let us state our main results in this paper. 
Let $\lambda\in P^+$ be given.  Let $(\mathscr{L}^{\lambda_+},\mathscr{B}^{\lambda_+})$ be a lower crystal base of $V_{m|0}(\lambda_+)$, and $(\mathscr{L}^{\lambda_-},\mathscr{B}^{\lambda_-})$ is a upper crystal base of $V_{0|n}(\lambda_-)$.
Set
\begin{equation*}\label{crystal base of Klambda}
\begin{split}
\mathscr{L}({K(\lambda)})&= \bigoplus_{S\subset \Phi^+_1}\A{\bf f}_S\otimes  \mathscr{L}^{\lambda_+}\otimes \mathscr{L}^{\lambda_-}  \subset K(\lambda),\\
\mathscr{B}({K(\lambda)})&=\left\{\ \pm\,{{\bf f}_S}\otimes b_+\otimes b_-\ \Big|\ S\subset \Phi^+_1,\ b_\pm\in\mathscr{B}^{\lambda^\pm}\ \right\} \subset  \mathscr{L}({K(\lambda)})/q\mathscr{L}({K(\lambda)}),
\end{split}
\end{equation*}
where we assume that $1_\lambda\in \mathscr{L}({K(\lambda)})_\lambda$.

\begin{thm}[\textsc{Existence}]\label{main result}
$\left(\mathscr{L}({K(\lambda)}),\mathscr{B}(K(\lambda))\right)$ is a crystal base of $K(\lambda)$. 
\end{thm}
\pf By \eqref{e0 action on PBW}, we have for $S\subset \Phi^+_1$
\begin{equation}\label{tilde e0 action on PBW}
\begin{split}
e_0'({\bf f}_S)&=
\begin{cases}
{\bf f}_{S\setminus\{\alpha_0\}}, & \text{if $\alpha_0\in S$}, \\
0, & \text{if $\alpha_0\not\in S$},
\end{cases}\ \  \ \
f_0({\bf f}_S)=
\begin{cases}
{\bf f}_{S\cup\{\alpha_0\}}, & \text{if $\alpha_0\not\in S$}, \\
0, & \text{if $\alpha_0\in S$}.
\end{cases}
\end{split}
\end{equation}
This implies that $\mathscr{L}({K(\lambda)})$ and $\mathscr{B}({K(\lambda)})\sqcup\{0\}$ are invariant under $\te_0$ and $\tf_0$, and $\tilde{f}_0 b = b'$ if and only if $b=\tilde{e}_0 b'$  for $b,b'\in \mathscr{B}(K(\lambda))$.
The other conditions in Definition \ref{def of crystal base of Kac module} follow from Propositions \ref{lower crystal base of K(0)}, \ref{upper crystal base of K(0)'} and \ref{lower = upper}. Hence $\left(\mathscr{L}({K(\lambda)}),\mathscr{B}(K(\lambda))\right)$ is a crystal base of $K(\lambda)$. 
\qed\vskip 3mm

Moreover, we have the following results.
\begin{thm}[\textsc{Connectedness}]\label{connectedness}
The crystal $\mathscr{B}(K(\lambda))/\{\pm 1\}$ is connected.
\end{thm}

\begin{cor}
We have
\begin{equation*}
\begin{split}
\mathscr{L}({K(\lambda)})&=\sum_{r\geq 0, k_1,\ldots,k_r\in I}\A \td{x}_{k_1}\cdots\td{x}_{k_r}1_\lambda, \\
\mathscr{B}(K(\lambda))&=\{\,\pm \,{\td{x}_{k_1}\cdots\td{x}_{k_r}1_\lambda}\!\!\! \mod{q\mathscr{L}({K(\lambda)})}\,|\,r\geq 0, k_1,\ldots,k_r\in I\,\}\setminus\{0\},
\end{split}
\end{equation*}
where $x=e$ or $f$ for each $k_i$.
\end{cor}

\begin{thm}[\textsc{Uniqueness}]\label{main result - uniqueness}
A crystal base of $K(\lambda)$ is unique up to a scalar multiplication.
\end{thm}

\begin{thm}[\textsc{Compatibility}]\label{main result - compatibility} Let $\lambda\in \td{P}^+$ be given. Let $\pi_\lambda : K(\lambda) \longrightarrow V(\lambda)$ be the $U$-module homomorphism such that $\pi_\lambda(1_\lambda)=v_\lambda$. 
Then  
\begin{itemize}
\item[(1)] $\pi_\lambda(\mathscr{L}(K(\lambda)))=\mathscr{L}(\lambda)$,

\item[(2)] $\ov{\pi}_\lambda(\mathscr{B}(K(\lambda)))= \mathscr{B}(\lambda) \cup\{0\}$, where $\ov{\pi}_\lambda : \mathscr{L}(K(\lambda))/q\mathscr{L}(K(\lambda)) \rightarrow \mathscr{L}(\lambda)/q\mathscr{L}(\lambda)$ is the induced $\mathbb{Q}$-linear map,

\item[(3)] $\ov{\pi}_\lambda$ restricts to a weight preserving bijection 
$$\ov{\pi}_\lambda : \{\,b\in \mathscr{B}(K(\lambda))\,|\,\ov{\pi}_\lambda(b)\neq 0 \,\}/\{\pm1\} \longrightarrow \mathscr{B}(\lambda)/\{\pm1\},$$ 
which commutes with $\te_k$ and $\tf_k$ for $k\in I$.
\end{itemize}
\end{thm}

The proof of Theorem \ref{connectedness} is given in Section \ref{Proof of connectedness} after a combinatorial description of $\mathscr{B}(K(\lambda))/\{\pm1\}$ in Section \ref{combinatorics of Kac module}. As in the case of an irreducible polynomial representation, $\mathscr{B}(K(\lambda))/\{\pm1\}$ may have a fake highest weight element, that is, there exists $b$ such that ${\rm wt}(b)\neq \lambda$ but $\te_k b= 0$ for all $k\in I$. Theorem \ref{main result - uniqueness} follows from  Theorem \ref{connectedness} and \cite[Lemma 2,7 (iii) and (iv)]{BKK}. The proof of Theorem \ref{main result - compatibility} is given in Section \ref{proof of compatibility}. 

\section{Combinatorics of crystals of Kac modules}
\subsection{Description of crystal operators on $\mathscr{B}({K(\lambda)})/\{\pm1\}$}\label{combinatorics of Kac module}
Fix $\lambda\in P^+$. The map sending $({\bf f}_S \otimes b_+\otimes b_-)$ to $(-S,  b_+, b_-)$ gives a bijection $$\mathscr{B}({K(\lambda)})/\{\pm1\}\ \longrightarrow\  \cP(\Phi^-_1)\times \mathscr{B}^{\lambda_+}\times \mathscr{B}^{\lambda_-}$$ as a set, where $\cP(\Phi^-_1)$ is the power set of $\Phi^-_1=-\Phi^+_1$, and $-S=\{-\beta\,|\,\beta\in S\}$ for $S\subset \Phi^+_1$. 

Note that $\cP(\Phi^-_1)$ can be identified with $\mathscr{B}({K(0)})/\{\pm1\}$ and its structure can be described as follows. Let $S\in \cP(\Phi^-_1)$ be given with $S=\{\,\beta_1\prec\ldots\prec\beta_r\,\}=\{\,\beta'_1\prec'\ldots\prec'\beta'_r\,\}$. Here, we assume that for $\alpha,\beta\in \Phi^-_1$, $\alpha\prec \beta$ (resp. $\alpha\prec' \beta$) if and only if $-\alpha\prec-\beta$ (resp. $-\alpha\prec'-\beta$).  For $k=0$, we have by \eqref{tilde e0 action on PBW}
\begin{equation*}\label{tilde e0  f0 action on K(0)}
\begin{split}
\te_0(S)&=
\begin{cases}
{S\setminus\{-\alpha_0\}}, & \text{if $-\alpha_0\in S$}, \\
0, & \text{if $-\alpha_0\not\in S$},
\end{cases} \ \  \ \
\tf_0(S)=
\begin{cases}
{S\cup\{-\alpha_0\}},  & \text{if $-\alpha_0\not\in S$}, \\
0, & \text{if $-\alpha_0\in S$}.
\end{cases}
\end{split}
\end{equation*}
Suppose that $k\neq 0$. By \eqref{sl2 action on basis of K}, we have for $i=1,\ldots,r$ 
\begin{equation*}
\begin{split}
\te_k(\beta_i)&=
\begin{cases}
\beta_i+\alpha_k, & \text{if $\beta_i+\alpha_k\in \Phi^-_1$}, \\
0, & \text{otherwise},
\end{cases} \ \ \ 
\tf_k(\beta_i)=
\begin{cases}
\beta_i-\alpha_k,  & \text{if $\beta_i-\alpha_k\in \Phi^-_1$}, \\
0, & \text{otherwise}.
\end{cases}
\end{split}
\end{equation*}
Then by \eqref{phi m|0} and \eqref{phi 0|n}, we have
\begin{equation*}
\td{x}_k S=
\begin{cases}
\{\,\gamma_1,\ldots, \gamma_r\,\}, & \text{if $k\in I_{m|0}$ and $\td{x}_k(\beta_1\otimes \cdots \otimes \beta_r)=\gamma_1\otimes \cdots \otimes \gamma_r$},\\
\{\,\gamma'_1,\ldots, \gamma'_r\,\}, & \text{if $k\in I_{0|n}$ and $\td{x}_k(\beta'_1\otimes \cdots \otimes \beta'_r)=\gamma'_1\otimes \cdots \otimes \gamma'_r$},\\
\end{cases}
\end{equation*}
($x=e, f$) following \eqref{lower tensor product rule} and \eqref{upper tensor product rule}. Here we assume that $\td{x}_k S=0$ if $\td{x}_k(\beta_1\otimes \cdots \otimes \beta_r)=0$ or $\td{x}_k(\beta'_1\otimes \cdots \otimes \beta'_r)=0$.

By Theorem \ref{main result}, we have the following.
\begin{prop}
For $(S,b_+,b_-)\in \mathscr{B}({K(\lambda)})/\{\pm1\}$ and $k\in I$,
\begin{equation*}
\tilde{x}_k(S,b_+,b_-)=
\begin{cases}
(S',b'_+,b_-), & \text{if $k\in I_{m|0}$ and $\tilde{x}_k(S
\otimes b_+ )=S'\otimes b'_+$,} \\
(S'',b_+,b_-''), & \text{if $k\in I_{0|n}$ and $\tilde{x}_k(S
\otimes b_- )=S''\otimes b''_-$},\\
(\tilde{x}_0S, b_+,b_-), & \text{if $k=0$},
\end{cases}
\end{equation*}
where $x=e,f$, and $\widetilde{x}_k(S,b_+,b_-)=0$ if any
of its components is $0$.
\end{prop}

Let $\mathscr{B}_+=\{\,\ov{m}<\ldots <\ov{1}\,\}$ and $\mathscr{B}_-=\{\,1<\ldots<n\,\}$, which are the subsets of even and odd elements in $\mathscr{B}$, respectively. Suppose that $\lambda\in P^+$ is given, where $\mu=(\lambda_{\ov{m}},\ldots,\lambda_{\ov{1}})$ and $\nu=(\lambda_{1},\ldots,\lambda_n)$ are partitions. By Theorem \ref{BKK}, we may identify $\mathscr{B}^{\lambda_+}$ with an $I_{m|0}$-colored subgraph $SST_{\mathscr{B}_+}(\mu)$ of $SST_{\mathscr{B}}(\mu)$ and $\mathscr{B}^{\lambda_-}$ with an $I_{0|n}$-colored subgraph $SST_{\mathscr{B}_-}(\nu')$ of $SST_{\mathscr{B}}(\nu')$.

\begin{ex}{\rm Suppose that $(m|n)=(3|3)$ and $\lambda=4\epsilon_{\ov{3}}+3\epsilon_{\ov{2}}+2\epsilon_{\ov{1}}+3\epsilon_1+\epsilon_2\in P^+$. Then  we may identify $\mathscr{B}^{\lambda_+}$ with $SST_{\mathscr{B}_+}(4,3,2)$ and $\mathscr{B}^{\lambda_-}$ with $SST_{\mathscr{B}_-}(2,1,1)$.
Consider the following triple $(S,U,V)\in \cP(\Phi^-_1)\times \mathscr{B}^{\lambda_+}\times \mathscr{B}^{\lambda_-}$, where \vskip 3mm
\begin{center}
$(S,U,V)=$$\Bigg(\{\,-\epsilon_{\ov{2}}+\epsilon_1 , -\epsilon_{\ov{2}}+\epsilon_2 , -\epsilon_{\ov{1}}+\epsilon_3\,\}$\ , \ 
${\def\lr#1{\multicolumn{1}{|@{\hspace{.6ex}}c@{\hspace{.6ex}}|}{\raisebox{-.3ex}{$#1$}}}\raisebox{-.6ex}
{$\begin{array}{cccc}
\cline{1-1}\cline{2-2}\cline{3-3}\cline{4-4}
\lr{\ov{3}}&\lr{\ov{3}}&\lr{\ov{3}}&\lr{\ov{2}}\\ 
\cline{1-1}\cline{2-2}\cline{3-3}\cline{4-4}
\lr{\ov{2}}&\lr{\ov{2}}& \lr{\ov{1}}\\ 
\cline{1-1}\cline{2-2} \cline{3-3}
\lr{\ov{1}}&\lr{\ov{1}}\\
\cline{1-1}\cline{2-2} 
\end{array}$}}$ \ , \ 
${\def\lr#1{\multicolumn{1}{|@{\hspace{.6ex}}c@{\hspace{.6ex}}|}{\raisebox{-.3ex}{$#1$}}}\raisebox{-.6ex}
{$\begin{array}{cccc}
\cline{1-1}\cline{2-2}
\lr{1} & \lr{3} \\
\cline{1-1} \cline{2-2}
\lr{2}\\
\cline{1-1} 
\lr{2}\\
\cline{1-1} 
\end{array}$}}\ \ \Bigg) $.
\end{center}\vskip 3mm

It is clear that $\te_0(S,U,V)=0$ and $\tf(S,U,V)=(S\cup\{-\alpha_0\},U,V)$. 

Let us compute $\tf_{\ov{2}}(S,U,V)$. Since $S=\{-\epsilon_{\ov{2}}+\epsilon_1\prec -\epsilon_{\ov{2}}+\epsilon_2\prec -\epsilon_{\ov{1}}+\epsilon_3\}$ and $\varphi_{\ov{2}}(S)=2>\varepsilon_{\ov{2}}(U)=1$, we have $\tf_{\ov{2}}(S\otimes U)=\left(\tf_{\ov{2}}S\right)\otimes U$ and hence \vskip 3mm
\begin{center}
$\tf_{\ov{2}}(S,U,V)=$$\Bigg(\{\,-\epsilon_{\ov{3}}+\epsilon_1 , -\epsilon_{\ov{2}}+\epsilon_2 , -\epsilon_{\ov{1}}+\epsilon_3\,\}$\ , \ 
${\def\lr#1{\multicolumn{1}{|@{\hspace{.6ex}}c@{\hspace{.6ex}}|}{\raisebox{-.3ex}{$#1$}}}\raisebox{-.6ex}
{$\begin{array}{cccc}
\cline{1-1}\cline{2-2}\cline{3-3}\cline{4-4}
\lr{\ov{3}}&\lr{\ov{3}}&\lr{\ov{3}}&\lr{\ov{2}}\\ 
\cline{1-1}\cline{2-2}\cline{3-3}\cline{4-4}
\lr{\ov{2}}&\lr{\ov{2}}& \lr{\ov{1}}\\ 
\cline{1-1}\cline{2-2} \cline{3-3}
\lr{\ov{1}}&\lr{\ov{1}}\\
\cline{1-1}\cline{2-2} 
\end{array}$}}$ \ , \ 
${\def\lr#1{\multicolumn{1}{|@{\hspace{.6ex}}c@{\hspace{.6ex}}|}{\raisebox{-.3ex}{$#1$}}}\raisebox{-.6ex}
{$\begin{array}{cccc}
\cline{1-1}\cline{2-2}
\lr{1} & \lr{3} \\
\cline{1-1} \cline{2-2}
\lr{2}\\
\cline{1-1} 
\lr{2}\\
\cline{1-1} 
\end{array}$}}\ \ \Bigg) $.
\end{center}\vskip 3mm

Next, let us compute $\tf_{2}(S,U,V)$. In this case, $S=\{-\epsilon_{\ov{1}}+\epsilon_3\prec'  -\epsilon_{\ov{2}}+\epsilon_2\prec' -\epsilon_{\ov{2}}+\epsilon_1\}$ and $\varphi_{2}(V)=1> \varepsilon_{2}(S)=0$, which implies that $\tf_2(S\otimes V)=S\otimes\left(\tf_2 V \right)$ (see \eqref{upper tensor product rule}). Therefore,\vskip 3mm
\begin{center}
$\tf_2(S,U,V)=$$\Bigg(\{\,-\epsilon_{\ov{2}}+\epsilon_1 , -\epsilon_{\ov{2}}+\epsilon_2 , -\epsilon_{\ov{1}}+\epsilon_3\,\}$\ , \ 
${\def\lr#1{\multicolumn{1}{|@{\hspace{.6ex}}c@{\hspace{.6ex}}|}{\raisebox{-.3ex}{$#1$}}}\raisebox{-.6ex}
{$\begin{array}{cccc}
\cline{1-1}\cline{2-2}\cline{3-3}\cline{4-4}
\lr{\ov{3}}&\lr{\ov{3}}&\lr{\ov{3}}&\lr{\ov{2}}\\ 
\cline{1-1}\cline{2-2}\cline{3-3}\cline{4-4}
\lr{\ov{2}}&\lr{\ov{2}}& \lr{\ov{1}}\\ 
\cline{1-1}\cline{2-2} \cline{3-3}
\lr{\ov{1}}&\lr{\ov{1}}\\
\cline{1-1}\cline{2-2} 
\end{array}$}}$ \ , \ 
${\def\lr#1{\multicolumn{1}{|@{\hspace{.6ex}}c@{\hspace{.6ex}}|}{\raisebox{-.3ex}{$#1$}}}\raisebox{-.6ex}
{$\begin{array}{cccc}
\cline{1-1}\cline{2-2}
\lr{1} & \lr{3} \\
\cline{1-1} \cline{2-2}
\lr{2}\\
\cline{1-1} 
\lr{3}\\
\cline{1-1} 
\end{array}$}}\ \ \Bigg) $.
\end{center}\vskip 3mm

}
\end{ex}

\begin{rem}\label{binary matrix}{\rm 
One may identify $\cP(\Phi^-_1)$ with 
$\M =\{\,A=(a_{ij})_{1\leq i\leq m, 1\leq j\leq n}\,|\,a_{ij}=0,1\,\}$.
The map sending $S$ to $A$ is a bijection from $\cP(\Phi^-_1)$ to $\M$, where $a_{ij}=1$ if and only if $-\epsilon_{\ov{i}}+\epsilon_j\in S$.
}
\end{rem}

\subsection{Connectedness of $\mathscr{B}(K(\lambda))/\{\pm1\}$}\label{Proof of connectedness}
Let $\mathscr{B}_+^\vee=\{\,\ov{1}^\vee,\ldots,\ov{m}^\vee\,\}$ be the crystal associated to the dual of the natural representation of $U_{m|0}$. We assume that $\mathscr{B}_+^\vee$ has a linear ordering $\ov{1}^\vee< \cdots <\ov{m}^\vee$. 
For a skew Young diagram $\mu/\nu$, we may define $\te_k$ and $\tf_k$ ($k\in I_{m|0}$) on $SST_{\mathscr{B}_+^\vee}(\mu/\nu)$ in the same way as in the case of $SST_{\mathscr{B}_+}(\mu/\nu)$ (cf.\cite{KashNaka}).

Let  $\lambda\in P^+$ be given. 
We assume that $\lambda_{\ov{m}}<0$ and $\lambda_{n}> 0$, and put 
\begin{equation*}
\mu=(\ell+\lambda_{\ov{m}},\ldots,\ell+\lambda_{\ov{1}}), \ \ \nu=(\lambda_1,\ldots,\lambda_n)'.
\end{equation*}
where $\ell$ is a positive integer such that $\ell+\lambda_{\ov{1}}>0$. We may identify $\mathscr{B}^{\lambda_+}$ with $SST_{\mathscr{B}^\vee_+}((\ell^m)/\mu)$, and $\mathscr{B}^{\lambda_-}$ with $SST_{\mathscr{B}_-}(\nu)$ so that $\mathscr{B}(K(\lambda))/\{\pm1\}$ can be identified with
$$ \cP(\Phi^-_1)\times SST_{\mathscr{B}^\vee_+}((\ell^m)/\mu)\times SST_{\mathscr{B}_-}(\nu).$$

Suppose that $S=\{\,\beta_1\prec\cdots\prec \beta_r\,\}\in \cP(\Phi^-_1)$ is given with $\beta_k=-\epsilon_{\ov{i_k}}+\epsilon_{j_k}$ for $1\leq k\leq r$. Let $w(S)=\ov{i_1}^\vee \ldots \ov{i_r}^\vee$, which is a finite word with alphabets in $\mathscr{B}_+^\vee$.   
For $U\in SST_{\mathscr{B}_+^\vee}((\ell^m)/\mu)$, we define
${ P}(U \leftarrow S)=(U \leftarrow w(S))$. Suppose that ${\rm sh}{P}(U\leftarrow
S)=(\ell^m)/\eta$ for some $\eta\subset \mu$. For $1\le k \le r$, let us
fill a box in $\mu/\eta$ with $j_k$ if it is created when
$\ov{i_k}^\vee$ is inserted into $(({S}\leftarrow{\ov{i_{r}}^\vee}\,)
\cdots)\leftarrow{\ov{i_{k+1}}^\vee}$. This defines the recording
tableau ${ Q}(U\leftarrow S)$. One can check that ${ Q}(U\leftarrow S)\in SST_{\mathscr{B}_+}(\mu/\eta)$ by modifying the proof of the dual RSK algorithm (cf.\cite{Ful,SS}), and that the correspondence from $(S,U)$ to $({ P}(U \leftarrow S),{ Q}(U \leftarrow S))$ is reversible.  Hence 
the map sending $(S,U,V)$ to $({P}(U \leftarrow S),{Q}(U \leftarrow S),V)$ gives a bijection
\begin{equation*}
\rho_\lambda : \mathscr{B}(K(\lambda))/\{\pm1\} \longrightarrow \mathscr{K}_\lambda,
\end{equation*}
where
\begin{equation*}
\mathscr{K}_\lambda=\bigsqcup_{\eta\subset \mu}SST_{\mathscr{B}_+^\vee}((\ell^m)/\eta)\times SST_{\mathscr{B}_-}( \mu/\eta)\times SST_{\mathscr{B}_-}(\nu).
\end{equation*}

Let us define $\te_k$ and $\tf_k$ on $\mathscr{K}_\lambda$ for $k\in I$.
The operators $\te_k$ and $\tf_k$ ($k\in I\setminus\{0\}$) are clearly defined on $\mathscr{K}_\lambda$ (for $k\in I_{0|n}$ we use the tensor product rule \eqref{upper tensor product rule} on the second and third components). Let us define $\te_0$ and $\tf_0$ on $\mathscr{K}_\lambda$. Let $(P,Q,V)\in\mathscr{K}_\lambda$ be given. For $k\geq 1$, let $a_k$ and $b_k$ be the top entries in the $k$th columns of $P$ and $Q$ (enumerated from the right), and let
\begin{equation*}
\sigma_k=
\begin{cases}
\ + & \text{if $a_k>\ov{1}^\vee$ or the $k$th column is empty},\\
\ - & \text{if $a_k=\ov{1}^\vee$ and $b_k=1$},\\
\ \ \cdot & \text{if otherwise}.\\
\end{cases}
\end{equation*}
Let $k_0$ be the smallest such that $\sigma_{k_0}\neq \ \cdot\ $. If $\sigma_{k_0}=+$, then we define $\te_0 (P,Q,V)=0$ and $\tf_0(P,Q,V)=(P',Q',V)$, where $(P',Q')$ is the pair of tableaux obtained from $(P,Q)$ by adding  $\boxed{\ov{1}^\vee}$ and $\boxed{1}$ on top of the $k_0$th columns of $P$ and $Q$, respectively. If $\sigma_{k_0}= -$, then we define $\tf_0 (P,Q,V)=0$ and $\te_0(P,Q,V)=(P',Q',V)$, where $(P',Q')$ is the pair of tableaux obtained from $(P,Q)$ by removing  $\boxed{\ov{1}^\vee}$ and $\boxed{1}$ in the $k_0$th columns of $P$ and $Q$, respectively. If $\sigma_{k}=\ \cdot \ $ for all $k$, then we define $\te_0(P,Q,V)=\tf_0(P,Q,V)=0$. Note that $\td{x}_0(P,Q,V)\in \mathscr{K}_\lambda$ if $\td{x}_0(P,Q,V)\neq 0$ ($x=e,f$).

\begin{lem}\label{connectedness special case}
$\rho_\lambda$ is a weight preserving bijection, which commutes with $\te_k$ and $\tf_k$ for $k\in I$.
\end{lem}
\pf By construction, $\rho_\lambda$ is a weight preserving bijection. So it remains to show that $\rho_\lambda$ commutes  with $\te_k$ and $\tf_k$ for $k\in I$. 

Let $(S,U,V)\in \mathscr{B}(K(\lambda))/\{\pm1\}$ be given with $\rho_\lambda(S,U,V)=(P,Q,V)$.

Suppose that $k\neq 0$. Then we claim that 
\begin{equation}\label{P commuting}
\td{x}_{i_1}\cdots \td{x}_{i_r}(S\otimes U)\neq 0 \ \ \Longleftrightarrow \ \ \td{x}_{i_1}\cdots \td{x}_{i_r}{P}\neq 0,
\end{equation}
for $i_1,\ldots,i_r\in I_{m|0}$ ($r\geq 1$), and
\begin{equation}\label{Q commuting}
\td{x}_{i_1}\cdots \td{x}_{i_r}S\neq 0\ \ \Longleftrightarrow \ \ \td{x}_{i_1}\cdots \td{x}_{i_r}{Q}\neq 0 ,
\end{equation}
for $i_1,\ldots,i_r\in I_{0|n}$ ($r\geq 1$), where $x=e, f$ for each $i_k$. 
We consider a binary matrix $M=(m_{ab})$ whose row indices are from $\mathscr{B}_+^\vee$ and column indices are from a linearly ordered set $\mathscr{B}_-\cup C$ for some linearly ordered  $\Z_2$-graded set $C$ with $|c|=1$ and $k<c$ for all $k\in\mathscr{B}_-$ and $c\in C$. Let $w(M)=(\ov{i_1}^\vee, j_1)\ldots (\ov{i_r}^\vee, j_r)$ be a biword such that
\begin{itemize}
\item[(1)] $\ov{i_p}^\vee\in \mathscr{B}^\vee_+$ and $j_p\in \mathscr{B}_-\cup C$ for $p=1,\ldots, r$,

\item[(2)] $j_{p}< j_{p+1}$ or ($i_p < i_{p+1}$ and $j_p=j_{p+1}$) for $p=1,\ldots, r-1$,

\item[(3)] $m_{ab}=1$ if and only if $(a,b)=(\ov{i_p}^\vee,j_p)$ for some $p=1,\ldots,r$.
\end{itemize}
Let ${\bf P}=(((\emptyset \leftarrow \ov{i_r}^\vee)\leftarrow \ov{i_{r-1}}^\vee) \cdots \leftarrow \ov{i_1}^\vee)$ and ${\bf Q}=(((\emptyset \leftarrow j_r)\leftarrow j_{r-1}) \cdots \leftarrow j_1)$, where $\emptyset$ is the empty tableau.
By the dual RSK correspondence, we have ${\bf P}\in SST_{\mathscr{B}_+^\vee}(\tau^\pi)$ and ${\bf Q}\in SST_{\mathscr{B}_-\cup C}(\tau^\pi)$ for some $\tau\in \mc{P}$.  
Suppose that 
\begin{itemize}
\item[(1)] $m_{ab}=1$ if and only if $-\epsilon_a+\epsilon_b \in S$ (see Remark \ref{binary matrix}),

\item[(2)] $((((\emptyset \leftarrow \ov{i_r}^\vee)\leftarrow \ov{i_{r-1}}^\vee) \cdots \leftarrow \ov{i_s}^\vee)=U$, where $j_p\in C$ if and only if $p=s,\ldots, r$. 
\end{itemize}
Then ${\bf P}=P$ and $Q$ is obtained by ignoring the entries from $C$ in the recording tableau ${\bf Q}$. 
Since the (dual) RSK correspondence is an isomorphism of bicrystals \cite{La} (see also Remark \ref{binary matrix} and \cite{K07}), we have \eqref{P commuting} and \eqref{Q commuting}. Therefore, $\rho_\lambda$ commutes  with $\te_k$ and $\tf_k$ for $k\in I\setminus\{0\}$.

Suppose that $k=0$.
We assume that $S=\{\,\beta_1\prec\cdots\prec \beta_r\,\}\in \mathscr{P}(\Phi^-_1)$ with $\beta_k=-\epsilon_{\ov{i_k}}+\epsilon_{j_k}$ for $1\leq k\leq r$, and hence $w(S)=\ov{i_1}^\vee \ldots \ov{i_r}^\vee$.

Suppose that $\te_0(P,Q,V)=(P',Q',V)\neq 0$. Let $k_0$ be such that $\sigma_k=\ \cdot \ $ for $k<k_0$ and   $\sigma_{k_0}=-$. So we have  $\boxed{\ov{1}^\vee}$ and $\boxed{1}$ in the $k_0$th columns of $P$ and $Q$, respectively.
Considering the bumping paths for each letter $\ov{i_k}^\vee$ in $U\leftarrow w(S)$, we see that $-\alpha_0\in S$ (that is, $\beta_1=-\epsilon_{\ov{1}}+\epsilon_1$ or $i_1=j_1=1$), and the insertion of $\ov{1}^\vee$ into  $(\cdots((U\leftarrow{\ov{i_r}^\vee})\leftarrow\ov{i_{r-1}}^\vee \,)\cdots \leftarrow \ov{i_2}^\vee)$ creates the pairs $\boxed{\ov{1}^\vee}$ and $\boxed{1}$ in the $k_0$th columns. This implies that 
\begin{equation*}
\begin{split}
P(U\leftarrow w(S\setminus\{-\alpha_0\}))=P',\ \
Q(U\leftarrow w(S\setminus\{-\alpha_0\}))=Q',
\end{split}
\end{equation*}
and
\begin{equation*}
\begin{split}
\rho_\lambda(\te_0(S,U,V))=\rho_\lambda(S\setminus\{-\alpha_0\},U,V)=(P',Q',V)=\te_0(P,Q,V).
\end{split}
\end{equation*}
Similarly, we can check that $\rho_\lambda(\tf_0(S,U,V))=\tf_0(P,Q,V)$ if $\tf_0(P,Q,V)\neq 0$.
Therefore, $\rho_\lambda$ commutes with $\te_0$ and $\tf_0$. \qed\vskip 3mm

\textsc{Proof of Theorem \ref{connectedness}}.
Let $D(k)$ denote the one dimensional $U$-module with weight $k\delta$ for $k\in \Z$. Since $K(\lambda)\otimes D(k)\cong K(\lambda+k\delta)$ for $k\in \Z$, $K(\lambda)\otimes D(k)$ has a crystal base, whose crystal can be identified with
\begin{equation*}
\cP(\Phi^-_1)\times \mathscr{B}^{\lambda_++k\delta_+}\times \mathscr{B}^{\lambda_-+k\delta_-}.
\end{equation*}
Recall that there exist bijections 
\begin{equation}\label{rational tableaux map}
\begin{split}
\sigma^k :  \mathscr{B}^{\lambda_+} \longrightarrow  \mathscr{B}^{\lambda_++k\delta_+},\ \ \ \
\tau^{k} :  \mathscr{B}^{\lambda_-} \longrightarrow  \mathscr{B}^{\lambda_-+k\delta_-},\\
\end{split}
\end{equation}
which commute with $\te_k$ and $\tf_k$ for $k\in I_{m|0}$ and $I_{0|n}$, respectively,  with ${\rm wt}(\sigma^k(b_+))={\rm wt}(b_+)+k\delta_+$ and ${\rm wt}(\tau^k(b_+))={\rm wt}(b_+)+k\delta_-$ for $b_{\pm}\in \mathscr{B}^{\lambda_{\pm}}$  (see \cite[Section 5.3]{K07} for more details). This impies that there exists a bijection
\begin{equation}\label{shift Klambda}
\varsigma^k : \mathscr{B}(K(\lambda))/\{\pm1\} \longrightarrow \cP(\Phi^-_1)\times \mathscr{B}^{\lambda_++k\delta_+}\times \mathscr{B}^{\lambda_-+k\delta_-},
\end{equation}
which commutes with $\te_k$ and $\tf_k$ ($k\in I$) and ${\rm wt}(\varsigma^k(b))={\rm wt}(b)+k\delta$.

Hence we may only consider the case when $\lambda_{\ov{m}}<0$ and $\lambda_{n}> 0$ (by taking $k\ll 0$ in \eqref{shift Klambda}).
By Lemma \ref{connectedness special case}, we may identify $\mathscr{B}(K(\lambda))/\{\pm1\}$ with $\mathscr{K}_{\lambda}$.

Now, we can apply the same argument in \cite[Theorem 4.8]{BKK} to prove that the crystal $\mathscr{B}(K(\lambda))/\{\pm1\}$ is connected 
(the only difference is that  a subtableau obtained from the first $m$ rows of $T\in \mathscr{B}(\lambda)/\{\pm1\}=SST_{\mathscr{B}}(\lambda^\circ)$  in the proof of \cite[Theorem 4.8]{BKK} is replaced with a pair of tableaux $(P,Q)$ in the first two components of $\mathscr{K}_\lambda$).
\qed

\subsection{Embedding of $\mathscr{B}(\lambda)/\{\pm1\}$ into $\mathscr{B}(K(\lambda))/\{\pm1\}$}\label{compatibility}
Let $\lambda\in\td{P}^+$ be given. By Theorem \ref{main result - compatibility}, there exists a unique injective map
\begin{equation}\label{embedding xi}
\xi_\lambda : \mathscr{B}(\lambda)/\{\pm1\} \longrightarrow \mathscr{B}(K(\lambda))/\{\pm1\}
\end{equation}
such that  for $b\in \mathscr{B}(\lambda)/\{\pm1\}$ and $k\in I$
\begin{itemize}
\item[(1)] ${\rm wt}(\xi_\lambda(b))={\rm wt}(b)$,

\item[(2)]   $\xi_\lambda(\td{x}_k b)=\td{x}_k\xi_\lambda (b)$ if $\td{x}_k b\neq 0$  $(x=e, f)$.
\end{itemize}
The purpose of  this subsection is to give an explicit description of $\xi_\lambda$.\vskip 2mm

Let  $\lambda^\circ=(\lambda^\circ_i)_{i\geq 1}$ be  the  partition in $\mc{P}_{m|n}$ corresponding to $\lambda$. 
We may identify $\mathscr{B}(\lambda)/\{\pm1\}$ with  $SST_{\mathscr{B}}(\lambda^\circ)$ by Theorem \ref{BKK}.  
Let 
\begin{equation*}
\begin{split}
\mu&= (\lambda^\circ_1,\ldots,\lambda^\circ_m), \ \
\nu=(\lambda^\circ_{m+1},\lambda^\circ_{m+2},\ldots).
\end{split}
\end{equation*}
For  $T\in SST_{\mathscr{B}}(\lambda^\circ)$, let 
\begin{itemize}
\item[$\cdot$] $T_{\leq m}$ \!\! : the subtableau of shape $\mu$ consisting of the first $m$ rows in $T$,

\item[$\cdot$] $T_{\leq m}^+$ : the subtableaux of $T_{\leq m}$ with entries in $\mathscr{B}_+$,

\item[$\cdot$] $T_{\leq m}^-$ : the subtableaux of $T_{\leq m}$ with entries in $\mathscr{B}_-$,

\item[$\cdot$] $T_{>m}$  : the complement of $T_{\leq m}$ in $T$.
\end{itemize}
Note that ${\rm sh}(T_{\leq m}^+)=\eta$ and ${\rm sh}(T_{\leq m}^-)=\mu/\eta$ for some $\eta=(\eta_1,\ldots,\eta_m)\in \mc{P}$, and  $T_{>m}\in SST_{\mathscr{B}_-}(\nu)$ can be regarded as an element in $\mathscr{B}^{\lambda_-}$.
\vskip 2mm

Let $\ell = \lambda^\circ_1=\mu_1$. Consider
\begin{equation*}
\sigma^{-\ell} : SST_{\mathscr{B}_+}(\eta) \longrightarrow SST_{\mathscr{B}_+^\vee}((\ell^m)/\eta)
\end{equation*}
which commutes with $\te_k$ and $\tf_k$ for $k\in I_{m|0}$ (see \eqref{rational tableaux map}).
Then the map sending $T$ to $(\sigma^{-\ell}(T_{\leq m}^+),T_{\leq m}^-,T_{> m})$ gives an injective map
\begin{equation}\label{embedding 1}
\imath_\lambda : \mathscr{B}(\lambda)/\{\pm1\} \longrightarrow \mathscr{K}_{\lambda -\ell\delta_+}.
\end{equation}
It is straightforward to check that (1) ${\rm wt}(\imath_\lambda(b))={\rm wt}(b)-\ell\delta_+$, (2)  $\imath_\lambda(\td{x}_k b)=\td{x}_k\imath_\lambda (b)$ if $\td{x}_k b\neq 0$  $(x=e, f)$ for $k\in I$.
 Since we have bijections
\begin{equation}\label{embedding 2}
\xymatrixcolsep{4pc}\xymatrixrowsep{4pc}\xymatrix{
\mathscr{K}_{\lambda -\ell\delta_+} \ar@{->}[r]^-{\rho_{\lambda-\ell\delta_+}^{-1}} 
& \mathscr{B}(K(\lambda-\ell\delta_+))/\{\pm1\}  \ar@{->}[r]^-{{\rm id}\times \sigma^\ell \times {\rm id}} & \mathscr{B}(K(\lambda))/\{\pm1\}},
\end{equation}
which commute with $\te_k$ and $\tf_k$ for $k\in I$,  by composing \eqref{embedding 1} and \eqref{embedding 2} we obtain a required map \eqref{embedding xi}.

\begin{ex}{\rm Suppose that $(m|n)=(3|3)$. Let $\lambda=4\epsilon_{\ov{3}}+3\epsilon_{\ov{2}}+2\epsilon_{\ov{1}}+2\epsilon_1$. Then $\lambda^\circ=(4,3,2,1,1)\in \mathcal{P}_{3|3}$. Consider\vskip 3mm

\begin{center}
$T=$ \ \ 
${\def\lr#1{\multicolumn{1}{|@{\hspace{.6ex}}c@{\hspace{.6ex}}|}{\raisebox{-.3ex}{$#1$}}}\raisebox{-.6ex}
{$\begin{array}{cccc}
\cline{1-1}\cline{2-2}\cline{3-3}\cline{4-4}
\lr{\ov{3}}&\lr{\ov{3}}&\lr{\ov{2}}&\lr{\ov{1}}\\ 
\cline{1-1}\cline{2-2}\cline{3-3}\cline{4-4}
\lr{\ov{2}}&\lr{\ov{1}}&\lr{3}\\ 
\cline{1-1}\cline{2-2}\cline{3-3}
\lr{1}&\lr{2} \\
\cline{1-1}\cline{2-2}
\lr{1}  \\
\cline{1-1} 
\lr{2}\\
\cline{1-1} 
\end{array}$}}$ \ \ $\in SST_{\mathscr{B}}(\lambda^\circ )$.
\end{center}\vskip 3mm
Keeping the above notations, we have\vskip 3mm
\begin{center}
$T^+_{\leq 3}=$ \ \ 
${\def\lr#1{\multicolumn{1}{|@{\hspace{.6ex}}c@{\hspace{.6ex}}|}{\raisebox{-.3ex}{$#1$}}}\raisebox{-.6ex}
{$\begin{array}{cccc}
\cline{1-1}\cline{2-2}\cline{3-3}\cline{4-4}
\lr{\ov{3}}&\lr{\ov{3}}&\lr{\ov{2}}&\lr{\ov{1}}\\ 
\cline{1-1}\cline{2-2}\cline{3-3}\cline{4-4}
\lr{\ov{2}}&\lr{\ov{1}}\\ 
\cline{1-1}\cline{2-2} 
\end{array}$}}$ \ \ , \ \ \ \
$T^-_{\leq 3}=$ \ \ 
${\def\lr#1{\multicolumn{1}{|@{\hspace{.6ex}}c@{\hspace{.6ex}}|}{\raisebox{-.3ex}{$#1$}}}\raisebox{-.6ex}
{$\begin{array}{cccc}
 \cline{3-3} 
& &\lr{3}\\ 
\cline{1-1}\cline{2-2}\cline{3-3}
\lr{1}&\lr{2} \\
\cline{1-1}\cline{2-2}
\end{array}$}}$ \ \ ,\ \ \ \ 
$T_{>3}=$ \ \ 
${\def\lr#1{\multicolumn{1}{|@{\hspace{.6ex}}c@{\hspace{.6ex}}|}{\raisebox{-.3ex}{$#1$}}}\raisebox{-.6ex}
{$\begin{array}{cccc}
\cline{1-1}\cline{2-2}
\lr{1}  \\
\cline{1-1} 
\lr{2}\\
\cline{1-1} 
\end{array}$}}$\ \ .
\end{center}\vskip 3mm
Then 
\begin{center}
$\sigma^{-4}\left(T^+_{\leq 3}\right)=$ \ \ 
${\def\lr#1{\multicolumn{1}{|@{\hspace{.6ex}}c@{\hspace{.6ex}}|}{\raisebox{-.3ex}{$#1$}}}\raisebox{-.6ex}
{$\begin{array}{cccc}
\cline{3-3}\cline{4-4}
 & &\lr{\ov{1}^\vee}&\lr{\ov{2}^\vee}\\ 
\cline{1-1}\cline{2-2}\cline{3-3}\cline{4-4}
\lr{\ov{1}^\vee}&\lr{\ov{2}^\vee}&\lr{\ov{3}^\vee}&\lr{\ov{3}^\vee}\\ 
\cline{1-1}\cline{2-2}\cline{3-3}\cline{4-4} 
\end{array}$}}$ \  
\end{center}\vskip 3mm
and $\imath_\lambda(T)=(\sigma^{-4}\left(T^+_{\leq 3}\right),T^-_{\leq 3},T_{>3})\in \mathscr{K}_{\lambda-4\delta_+}$. Applying $\rho^{-1}_{\lambda-4\delta_+}$ to this triple (see the proof of Lemma \ref{connectedness special case}), we get $(S,U,T_{>3})\in \mathscr{B}(K(\lambda-4\delta_+))/\{\pm1\}$, where\vskip 3mm

\begin{center}
$S=\{\,-\epsilon_{\ov{3}}+\epsilon_1 , -\epsilon_{\ov{2}}+\epsilon_2 , -\epsilon_{\ov{1}}+\epsilon_3\,\}$ \ ,\ \
$U=$ \ \ \
${\def\lr#1{\multicolumn{1}{|@{\hspace{.6ex}}c@{\hspace{.6ex}}|}{\raisebox{-.3ex}{$#1$}}}\raisebox{-.6ex}
{$\begin{array}{ccc}
\cline{2-2} 
 &  \lr{\ov{1}^\vee}\\ 
\cline{1-1}\cline{2-2} 
 \lr{\ov{2}^\vee}&\lr{\ov{3}^\vee} \\ 
\cline{1-1}\cline{2-2} 
\end{array}$}}$ \  .
\end{center}\vskip 3mm
Applying $\sigma^4$ to $U$, we have  \vskip 3mm
\begin{center}
$\xi_\lambda(T)=$ $\Bigg(\{\,-\epsilon_{\ov{3}}+\epsilon_1 , -\epsilon_{\ov{2}}+\epsilon_2 , -\epsilon_{\ov{1}}+\epsilon_3\,\}$\ , \ 
${\def\lr#1{\multicolumn{1}{|@{\hspace{.6ex}}c@{\hspace{.6ex}}|}{\raisebox{-.3ex}{$#1$}}}\raisebox{-.6ex}
{$\begin{array}{cccc}
\cline{1-1}\cline{2-2}\cline{3-3}\cline{4-4}
\lr{\ov{3}}&\lr{\ov{3}}&\lr{\ov{3}}&\lr{\ov{2}}\\ 
\cline{1-1}\cline{2-2}\cline{3-3}\cline{4-4}
\lr{\ov{2}}&\lr{\ov{2}}& \lr{\ov{1}}\\ 
\cline{1-1}\cline{2-2} \cline{3-3}
\lr{\ov{1}}&\lr{\ov{1}}\\
\cline{1-1}\cline{2-2} 
\end{array}$}}$ \ , \ 
${\def\lr#1{\multicolumn{1}{|@{\hspace{.6ex}}c@{\hspace{.6ex}}|}{\raisebox{-.3ex}{$#1$}}}\raisebox{-.6ex}
{$\begin{array}{cccc}
\cline{1-1}\cline{2-2}
\lr{1}  \\
\cline{1-1} 
\lr{2}\\
\cline{1-1} 
\end{array}$}}\ \ \Bigg) $,
\end{center}
which belongs to $\mathscr{B}(K(\lambda))/\{\pm1\}$.
}
\end{ex}

\section{Compatibility with crystals of polynomial representations}\label{proof of compatibility}
In this section, we give a proof of Theorem \ref{main result -  compatibility}. We fix $\lambda\in \td{P}^+$ throughout this section. 

\subsection{}\label{step 1} Since $K(\lambda)$ is completely reducible as a $U_{m,n}$-module, we have  
\begin{equation*}
K(\lambda)=\bigoplus_{i\in X_\lambda}K(\lambda;\sigma_i),
\end{equation*}
where $X_\lambda$ is an index set and $K(\lambda;\sigma_i)$ is an irreducible $U_{m,n}$-submodule with highest weight $\sigma_i$. Let $X'_\lambda\subset X_\lambda$ be such that ${\rm Ker}\, \pi_\lambda =\bigoplus_{i\in X'_\lambda}K(\lambda;\sigma_i)$, and put $Y_\lambda=X_\lambda\setminus X'_\lambda$. We have as a $U_{m,n}$-module
\begin{equation*}
V(\lambda)=\bigoplus_{i\in Y_\lambda}V(\lambda;\sigma_i), 
\end{equation*}
where  $V(\lambda;\sigma_i)= \pi_\lambda(K(\lambda;\sigma_i))\cong K(\lambda;\sigma_i)$ as a $U_{m,n}$-module for $i\in Y_\lambda$. 
For $i\in X_\lambda$, put  
\begin{equation*}
\begin{split}
\mathscr{L}(K(\lambda;\sigma_i))&=\mathscr{L}(K(\lambda))\cap K(\lambda;\sigma_i), \\ 
\mathscr{B}(K(\lambda;\sigma_i))&=\mathscr{B}(K(\lambda))\cap \left(\mathscr{L}(K(\lambda;\sigma_i))/q\mathscr{L}(K(\lambda;\sigma_i)) \right).
\end{split}
\end{equation*}
By Propositions \ref{lower crystal base of K(0)}, \ref{upper crystal base of K(0)'} and \ref{lower = upper}, $(\mathscr{L}(K(\lambda;\sigma_i)),\mathscr{B}(K(\lambda;\sigma_i))/\{\pm1 \})$ is  a lower crystal base of $K(\lambda;\sigma_i)$ as a $U_{m|0}$-module and a upper crystal base as a $U_{0|n}$-module, and 
\begin{equation*}\label{U0 decomp 1}
\begin{split}
\mathscr{L}(K(\lambda)) &=\bigoplus_{i\in X_\lambda}\mathscr{L}(K(\lambda;\sigma_i)), \ \ 
\mathscr{B}(K(\lambda)) =\bigsqcup_{i\in X_\lambda}\mathscr{B}(K(\lambda;\sigma_i)),
\end{split}
\end{equation*}
(see \cite[Lemma 2.6.3]{Kash1}). 
Let $$\mathscr{L}(\lambda)'=\pi_\lambda(\mathscr{L}(K(\lambda))),\ \ \mathscr{B}(\lambda)'=\ov{\pi}_\lambda(\mathscr{B}(K(\lambda))),$$ where $\ov{\pi}_\lambda : \mathscr{L}(K(\lambda))/q\mathscr{L}(K(\lambda)) \rightarrow \mathscr{L}(\lambda)'/q \mathscr{L}(\lambda)'$ is the induced $\mathbb{Q}$-linear map. For $i\in Y_\lambda$, put
\begin{equation*}
\begin{split}
\mathscr{L}(\lambda;\sigma_i)'&=\pi_\lambda(\mathscr{L}(K(\lambda;\sigma_i))), \ \
\mathscr{B}(\lambda;\sigma_i)'=\ov{\pi}_\lambda(\mathscr{B}(K(\lambda;\sigma_i))).
\end{split}
\end{equation*}
Then $(\mathscr{L}(\lambda;\sigma_i)',\mathscr{B}(\lambda;\sigma_i)'/\{\pm1\})$ is a s  a lower crystal base of $V(\lambda;\sigma_i)$ as a $U_{m|0}$-module and a upper crystal base as a $U_{0|n}$-module, and 
\begin{equation}\label{U0 decomp 2}
\mathscr{L}(\lambda)'=\bigoplus_{i\in Y_\lambda}\mathscr{L}(\lambda;\sigma_i)', \ \ \mathscr{B}(\lambda)'=\bigsqcup_{i\in Y_\lambda}\mathscr{B}(\lambda;\sigma_i)'.
\end{equation}


\subsection{}\label{step 2} For $\ell\in \Z$, let $$\theta_\ell : K(\lambda) \longrightarrow K(\lambda+\ell\delta_+)$$ be a $U^-$-linear map such that $\theta_{\ell}(1_\lambda)=1_{\lambda+\ell\delta_+}$. It is well-defined and indeed an isomorphism of $U^-$-modules. If we identify $V_{m|0}(\lambda_++\ell\delta_+)$ with $V_{m|0}(\lambda_+)\otimes D^\ell$, where $D^\ell$ is the one dimensional $U_{m|0}$-module with highest weight vector $v_{\ell\delta_+}$ of weight $\ell\delta_+$, then we have 
\begin{equation*}
\theta_{\ell}(u\otimes v_+\otimes v_-)=u\otimes (v_+\otimes v_{\ell\delta_+})\otimes v_-
\end{equation*}
for $u\in K$, $v_+\in V_{m|0}(\lambda_+)$ and $v_-\in V_{0|n}(\lambda_-)$.  In particular, we have $\theta_{\ell}(\mathscr{L}(K(\lambda)))=\mathscr{L}(K(\lambda+\ell\delta_+))$.
Let $$S(\lambda,\ell\delta_+) : V(\lambda)\otimes V(\ell\delta_+) \longrightarrow V(\lambda)$$ be a $U^-$-linear map given by
$S(\lambda,\ell\delta_+)(v\otimes v_{\ell\delta_+})=v$ and $S(\lambda,\ell\delta_+)(v\otimes w)=0$ 
for $v\in V(\lambda)$ and $w\in V(\ell\delta_+)\setminus V(\ell\delta_+)_{\ell\delta_+}$. 
Let $$\Phi(\lambda,\ell\delta_+) : V(\lambda+\ell\delta_+) \longrightarrow V(\lambda)\otimes V(\ell\delta_+)$$ be an injective $U$-module homomorphism given by $\Phi(\lambda,\ell\delta_+)(v_{\lambda+\ell\delta_+})=v_\lambda\otimes v_{\ell\delta_+}$. Since $V(\lambda)\otimes V(\ell\delta_+)$ is completely reducible \cite[Theorem 2.12]{BKK} and $U(v_\lambda\otimes v_{\ell\delta_+})\cong V(\lambda+\ell\delta_+)$, $\Phi(\lambda,\ell\delta_+)$ is well-defined. 
We put $$\vartheta_{-\ell} = S(\lambda,\ell\delta_+)\circ\Phi(\lambda,\ell\delta_+) : V(\lambda+\ell\delta_+) \longrightarrow V(\lambda),$$ which is $U^-$-linear.
Now, we have the following commutative diagram
\begin{equation*}\label{commuting diagram}
\xymatrixcolsep{4pc}\xymatrixrowsep{4pc}\xymatrix{
K(\lambda+\ell\delta_+) \ar@{->}[d]^{\theta_{-\ell}} \ar@{->}[r]^{\pi_{\lambda+\ell\delta_+}} 
& V(\lambda+\ell\delta_+) \ar@{->}[d]^{\vartheta_{-\ell}} \ar@{->}[r]^-{\Phi(\lambda,\ell\delta_+)} & V(\lambda)\otimes V(\ell\delta_+)\ar@{->}[ld]^{S(\lambda,\ell\delta_+)}\\
K(\lambda) \ar@{->}[r]^{\pi_\lambda} 
& V(\lambda)   }
\end{equation*}
where the vertical maps are $U^-$-linear and the horizontal ones are $U$-linear. Note that if we take $\ell\gg 0$ such that  the partition  corresponding to $\lambda+\ell\delta_+$ contains the rectangle $(n^m)$, then ${\rm ch}K(\lambda+\ell\delta_+)={\rm ch}V(\lambda+\ell\delta_+)$ \cite[Theorm 6.20]{BR} and hence $K(\lambda+\ell\delta_+)=V(\lambda+\ell\delta_+)$ or $\pi_{\lambda+\ell\delta_+}$ is an isomorphism.

\subsection{}\label{step 3} We will first prove   Theorem \ref{main result - compatibility} for $\lambda+\ell\delta_+\in \td{P}^+$ with $\ell\gg 0$. 
For $b\in \mathscr{B}(K(\lambda))/\{\pm1\}$, let  $u(b)\in U^-$ be a homogeneous element such that $u(b)1_\lambda\in \mathscr{L}(K(\lambda))$ and $u(b)1_\lambda\equiv b \mod{q\mathscr{L}(K(\lambda))}$. By Nakayama's lemma, $\{\,u(b)1_\lambda\,|\,b\in \mathscr{B}(K(\lambda))/\{\pm1\}\,\}$ is an $\A$-basis of $\mathscr{L}(K(\lambda))$ (and hence a $\mathbb{Q}(q)$-basis of $K(\lambda)$). Moreover, $\{\,\theta_{\ell}(u(b)1_\lambda)=u(b)1_{\lambda+\ell\delta_+}\,|\,b\in\mathscr{B}(K(\lambda))/\{\pm1\}\,\}$ is an $\A$-basis of $\mathscr{L}(K(\lambda+\ell\delta_+))$ since $\theta_{\ell}$ is a $U^-$-linear isomorphism with $\theta_{\ell}(\mathscr{L}(K(\lambda)))=\mathscr{L}(K(\lambda+\ell\delta_+))$. So, we can take a set of homogeneous vectors
\begin{equation*}\label{U lambda}
\mathscr{U}(\lambda)=\{\,u_1,\ldots, u_d\,\}\subset U^-
\end{equation*}
such that $\{\,u_1\,1_{\lambda+\ell\delta_+},\ldots, u_d\,1_{\lambda+\ell\delta_+}\,\}$ is an $\A$-basis of $\mathscr{L}(K(\lambda+\ell\delta_+))$ for $\ell\in\Z$, where $d=\dim K(\lambda)$.

Now, we choose $\ell\gg 0$ such that
\begin{itemize}
\item[(1)] $V(\lambda+\ell\delta_+)=K(\lambda+\ell\delta_+)$,

\item[(2)] $q^{2\ell}\,t_0^2\,e''_0(u)1_{\lambda}\in q\mathscr{L}(K(\lambda))$ for all $u\in \mathscr{U}(\lambda)$.
\end{itemize}
Since $t_0^2\,e''_0(u)1_{\lambda}=q^{2\langle h_0, \lambda+\alpha+\alpha_0\rangle}e''_0(u)1_{\lambda}$ with ${\rm wt}(u)=\alpha$, the condition (2) implies
\begin{equation}\label{e'' condition}
\begin{split}
t_0^2\,e''_0(u)1_{\lambda+\ell\delta_+}&= q^{2\langle h_0,\lambda+\alpha+\alpha_0 +\ell\delta_+ \rangle}\,e''_0(u)1_{\lambda+\ell\delta_-} \\  
&=q^{2\langle h_0,\lambda+\alpha+\alpha_0  \rangle +2\ell}\,e''_0(u)1_{\lambda+\ell\delta_+} \\
 &=\theta_{\ell}(q^{2\ell}\,t_0^2\,e''_0(u)1_{\lambda}) \in q\mathscr{L}(K(\lambda+\ell\delta_+)).
\end{split}
\end{equation}
We claim that 
\begin{equation}\label{invariance for large k}
\te_k \mathscr{L}(\lambda+\ell\delta_+)'\subset \mathscr{L}(\lambda+\ell\delta_+)',\ \ \ \tf_k \mathscr{L}(\lambda+\ell\delta_+)'\subset \mathscr{L}(\lambda+\ell\delta_+)'
\end{equation}
for $k\in I$. By (1), we have $\mathscr{L}(\lambda+\ell\delta_+)'=\pi_{\lambda+\ell\delta_+}(\mathscr{L}(K(\lambda+\ell\delta_+)))=\mathscr{L}(K(\lambda+\ell\delta_+))$, but to emphasize that the crystal operators in \eqref{invariance for large k} are those on the modules in $\mathcal{O}_{int}$  (that is, $\te_0=q^{-1}t_0e_0$  in \eqref{e0 f0 for Oint}), we use the notation $\mathscr{L}(\lambda+\ell\delta_+)'$.  

It is clear that $\mathscr{L}(\lambda+\ell\delta_+)'$ is invariant under $\te_k$ and $\tf_k$ for $k\in I\setminus \{0\}$ since $\mathscr{L}(\lambda+\ell\delta_+)'$ is a lower crystal lattice as a $U_{m|0}$-module and a upper crystal lattice as a $U_{0|n}$-module  (see \eqref{U0 decomp 2}). Also we have $\tf_0\mathscr{L}(\lambda+\ell\delta_+)'=f_0\mathscr{L}(\lambda+\ell\delta_+)'\subset \mathscr{L}(\lambda+\ell\delta_+)'$ since $\mathscr{L}(K(\lambda+\ell\delta_+))$ is invariant under $\tf_0=f_0$. So it remains to show that $\mathscr{L}(\lambda+\ell\delta_+)'$ is invariant under $\te_0$.  Let $u\in \mathscr{U}(\lambda)$ be given with ${\rm wt}(u)=\alpha$. Then 
\begin{equation*}
e_0 u = (-1)^{|\alpha|} u e_0 +\frac{t_0e''_0(u)-t_0^{-1}e'_0(u)}{q-q^{-1}}
\end{equation*}
or
\begin{equation*}
q^{-1}t_0e_0 u = (-1)^{|\alpha|} q^{-1}t_0 u e_0 +\frac{e'_0(u)-t^2_0e''_0(u)}{1-q^2},
\end{equation*}
which implies that
\begin{equation}\label{compatibility of e0}
\te_0(u 1_{\lambda+\ell\delta_+}) =q^{-1}t_0e_0 u 1_{\lambda+\ell\delta_+}= \frac{1}{1- q^2}\left(e'_0(u) 1_{\lambda+\ell\delta_+}- t^2_0e''_0(u) 1_{\lambda+\ell\delta_+}\right).
\end{equation}
Since $u \,1_{\lambda+\ell\delta_+}\in \mathscr{L}(K(\lambda+\ell\delta_+))$, we have  $ e'_0(u) 1_{\lambda+\ell\delta_+} \in \mathscr{L}(K(\lambda+\ell\delta_+))$ by \eqref{e0 f0 for K}. By \eqref{e'' condition}, we have $ t^2_0e''_0(u) 1_{\lambda+\ell\delta_+}\in q\mathscr{L}(K(\lambda+\ell\delta_+))$. Therefore, $q^{-1}t_0e_0 u 1_{\lambda+\ell\delta_+} \in \mathscr{L}(K(\lambda+\ell\delta_+))$. This proves that 
$\te_0 \mathscr{L}(\lambda+\ell\delta_+)'\subset \mathscr{L}(\lambda+\ell\delta_+)'$.

By \eqref{invariance for large k}, $\mathscr{L}(\lambda+\ell\delta_+)'$ is a crystal lattice of $V(\lambda+\ell\delta_+)$, and by the uniqueness of a crystal lattice \cite[Lemma 2.7 (iii)]{BKK}, we have 
\begin{equation}\label{condition (1) for special}
\mathscr{L}(\lambda+\ell\delta_+)'=\mathscr{L}(\lambda+\ell\delta_+).
\end{equation}
It is clear that  the induced map $\ov{\pi}_{\lambda+\ell\delta_+} : \mathscr{L}(K(\lambda+\ell\delta_+))/q\mathscr{L}(K(\lambda+\ell\delta_+)) \rightarrow \mathscr{L}(\lambda+\ell\delta_+)/q\mathscr{L}(\lambda+\ell\delta_+)$ commutes with $\te_k$, $\tf_k$ and $\tf_0$ ($k\in I\setminus \{0\}$). Also by \eqref{compatibility of e0},  we have $\te_0(u 1_{\lambda+\ell\delta_+})\equiv e'_0(u) 1_{\lambda+\ell\delta_+}$$\mod{q\mathscr{L}(\lambda+\ell\delta_+)}$ and hence $\ov{\pi}_{\lambda+\ell\delta_+}$ commutes with $\te_0$. Therefore, $\ov{\pi}_{\lambda+\ell\delta_+}(\mathscr{B}(K(\lambda+\ell\delta_+)))=\mathscr{B}(\lambda+\ell\delta_+)$ by Theorem \ref{connectedness}, and we have a weight preserving bijection 
\begin{equation}\label{condition (2) for special}
\ov{\pi}_{\lambda+\ell\delta_+} : \mathscr{B}(K(\lambda+\ell\delta_+))/\{\pm1\} \longrightarrow \mathscr{B}(\lambda+\ell\delta_+)/\{\pm1\},
\end{equation}
which commutes with $\te_k$ and $\tf_k$ for $k\in I$. By \eqref{condition (1) for special} and \eqref{condition (2) for special}, Theorem \ref{main result - compatibility} holds for $\lambda+\ell\delta_+$ for $\ell\gg 0$.
 
\subsection{}\label{step 4} Let $\ell$ be as in Section \ref{step 3}. Since $\Phi(\lambda,\ell\delta_+)(\mathscr{L}(\lambda+\ell\delta_+))\subset \mathscr{L}(\lambda)\otimes \mathscr{L}(\ell\delta_+)$ and $S(\lambda,\ell\delta_+)(\mathscr{L}(\lambda)\otimes \mathscr{L}(\ell\delta_+))=\mathscr{L}(\lambda)$, we have the following commutative diagram
\begin{equation*}
\xymatrixcolsep{4pc}\xymatrixrowsep{4pc}\xymatrix{
\mathscr{L}(\lambda+\ell\delta_+)/q\mathscr{L}(\lambda+\ell\delta_+) \ar@{->}[d]^{\ov{\vartheta_{-\ell}}} \ar@{->}[r]^-{\ov{\Phi(\lambda,\ell\delta_+)}} & \mathscr{L}(\lambda)\otimes \mathscr{L}(\ell\delta_+)/q\mathscr{L}(\lambda)\otimes \mathscr{L}(\ell\delta_+)\ar@{->}[ld]^{\ov{S(\lambda,\ell\delta_+)}}\\
\mathscr{L}(\lambda)/q\mathscr{L}(\lambda) }
\end{equation*}
Since 
\begin{equation*}
\begin{split}
&\ov{\Phi(\lambda,\ell\delta_+)}(\mathscr{B}(\lambda+\ell\delta_+)/\{\pm1\})\subset \mathscr{B}(\lambda)\otimes \mathscr{B}(\ell\delta_+)/\{\pm1\},\\
&\ov{S(\lambda,\ell\delta_+)}(\mathscr{B}(\lambda)\otimes \mathscr{B}(\ell\delta_+)/\{\pm1\})\subset \mathscr{B}(\lambda)/\{\pm1\}\cup\{0\},
\end{split}
\end{equation*}
$\ov{\vartheta_{-\ell}}$ induces a map 
\begin{equation*}
\ov{\vartheta_{-\ell}} : \mathscr{B}(\lambda+\ell\delta_+)/\{\pm 1\} \longrightarrow \mathscr{B}(\lambda)/\{\pm 1\}\cup\{0\}.
\end{equation*}
Note that $\ov{\Phi(\lambda,\ell\delta_+)}$ is an injecttive map which commutes with $\te_k$ and $\tf_k$ for $k\in I$ up to a multiplication by $\pm1 $.

Let us describe $\ov{\vartheta_{-\ell}}$ more explicitly. Let $\lambda^\circ$ and $(\lambda+\ell\delta_+)^\circ$ be the partitions in $\mathcal{P}_{m|n}$ corresponding to $\lambda$ and $\lambda+\ell\delta_+$, respectively. 
Let us identify $\mathscr{B}(\lambda)/\{\pm 1\}$ and $\mathscr{B}(\lambda+\ell\delta_+)/\{\pm 1\}$ with $SST_{\mathscr{B}}(\lambda^\circ)$ and $SST_{\mathscr{B}}((\lambda+\ell\delta_+)^\circ)$, respectively. Also, we may identify $\mathscr{B}(\ell\delta_+)/\{\pm 1\}$ with $SST_{\mathscr{B}}(\ell^m)$.

Suppose that $T\in SST_{\mathscr{B}}((\lambda+\ell\delta_+)^\circ)$ is given. By \cite[Theorem 4.18]{KK}, the multiplicity of $V(\lambda+\ell\delta_+)$ in $V(\lambda)\otimes V(\ell\delta_+)$ is 1, and there exist unique $T_1\in SST_{\mathscr{B}}(\lambda^\circ)$ and $T_2\in SST_{\mathscr{B}}(\ell^m)$ such that $T_1\otimes T_2$ generates the same $I$-colored oriented graph as that of $T$ (called crystal equivalent in \cite{BKK}). Indeed, we have $T=(T_2\rightarrow T_1):=(\psi(T_2)\rightarrow T_1)$, where $\psi$ is as in \eqref{admissible reading}.  This implies that $\ov{\Phi(\lambda,\ell\delta_+)}(T)=T_1\otimes T_2$. By definition of $S(\lambda,\ell\delta_+)$, we have
\begin{equation*}
\ov{S(\lambda,\ell\delta_+)}(T_1\otimes T_2)=
\begin{cases}
T_1, & \text{if $T_2 = H_{\ell\delta_+}$}, \\
0, & \text{otherwise}, \\
\end{cases}
\end{equation*}
where $H_{\ell\delta_+}\in SST_{\mathscr{B}}(\ell^m)$ is the highest weight element with weight $\ell\delta_+$ with $H_{\ell\delta_+}(i,j)=\ov{m-i+1}$ for $1\leq i\leq m$ and $1\leq j\leq \ell$.  Hence, we have
\begin{equation*}
\ov{\vartheta_{-\ell}}(T)=
\begin{cases}
T_1, & \text{if $T_2 = H_{\ell\delta_+}$}, \\
0, & \text{otherwise}. \\
\end{cases}
\end{equation*}
In particular, $\ov{\vartheta_{-\ell}}$ is surjective. By Nakayama's lemma, $\vartheta_{-\ell}(\mathscr{L}(\lambda+\ell\delta_+))=\mathscr{L}(\lambda)$, and
\begin{equation*}
\begin{split}
\pi_{\lambda}(\mathscr{L}(K(\lambda)))&=\vartheta_{-\ell}\circ \pi_{\lambda+\ell\delta_+}\circ\theta_{\ell}(\mathscr{L}(K(\lambda))) \\
&=\vartheta_{-\ell}\circ \pi_{\lambda+\ell\delta_+}(\mathscr{L}(K(\lambda+\ell\delta_+)))\\
&=\vartheta_{-\ell}(\mathscr{L}(\lambda+\ell\delta_+))=\mathscr{L}(\lambda).\\
\end{split}
\end{equation*}
This proves  Theorem \ref{main result - compatibility} (1) and (2).

Another way to describe $\ov{\vartheta_{-\ell}}(T)$ is as follows; for $T\in SST_{\mathscr{B}}((\lambda+\ell\delta_+)^\circ)$, let $T_1\in SST_{\mathscr{B}}(\lambda^\circ)$ and $T_2\in SST_{\mathscr{B}}(\ell^m)$ be the unique tableaux such that $(T_2\rightarrow T_1)=T$. Let $T'$ be the sub tableau in $SST_{\mathscr{B}}(\ell^m)$ with $T'(i,j)=T(i,j)$ for $1\leq i\leq m$ and $1\leq j\leq \ell$. By considering  the recording tableau of $(T_2\rightarrow T_1)$, we see that
$T'=H_{\ell\delta_+}$ if and only if $T_2=H_{\ell\delta_+}$, and in this case $T_1$ is given by $T_1(i,j)=T(i,j+\ell)$ for $1\leq i\leq m$ and $1\leq j\leq \lambda^\circ_i$, and $T_1(i,j)=T(i,j)$ for $i\geq m$ and $1\leq j\leq \lambda^\circ_i$.  

Let $T\in SST_{\mathscr{B}}((\lambda+\ell\delta_+)^\circ)$ be given such that $\ov{\vartheta_{-\ell}}(T)\neq 0$ and $\ov{\vartheta_{-\ell}}(\td{x}_kT)\neq 0$ for some $k\in I$ ($x=e,f$).
Suppose that $\ov{\Phi(\lambda,\ell\delta_+)}(T)=T_1\otimes H_{\ell\delta_+}$ for some $T_1$.  Since $\ov{\vartheta_{-\ell}}(T)\neq 0$, $\ov{\vartheta_{-\ell}}(\td{x}_kT)\neq 0$ and $\ov{\Phi(\lambda,\ell\delta_+)}(\td{x}_kT)=\td{x}_k\ov{\Phi(\lambda,\ell\delta_+)}(T)$, we have $\td{x}_k(T_1\otimes H_{\ell\delta_+})=(\td{x}_kT_1)\otimes H_{\ell\delta_+}$, and hence 
\begin{equation}\label{commuting property of vartheta}
\ov{\vartheta_{-\ell}}(\td{x}_kT)=\td{x}_k\ov{\vartheta_{-\ell}}(T).
\end{equation}

Finally, we have a commutative diagram
\begin{equation*}
\xymatrixcolsep{4pc}\xymatrixrowsep{4pc}\xymatrix{
\mathscr{B}(K(\lambda+\ell\delta_+))/\{\pm1\} \ar@{->}[d]^{\ov{\theta_{-\ell}}} \ar@{->}[r]^-{\ov{\pi_{\lambda+\ell\delta_+}}} 
& \mathscr{B}(\lambda+\ell\delta_+)/\{\pm1\} \ar@{->}[d]^{\ov{\vartheta_{-\ell}}} \\
\mathscr{B}(K(\lambda))/\{\pm1\} \ar@{->}[r]^{\ov{\pi_\lambda}} 
& \mathscr{B}(\lambda)/\{\pm1\}   }
\end{equation*}
Here  $\ov{\theta_{-\ell}}$ and $\ov{\pi_{\lambda+\ell\delta_+}}$ are bijections commuting with $\te_k$ and $\tf_k$ for $k\in I$. By \eqref{commuting property of vartheta}, we conclude that Theorem \ref{main result - compatibility} (3) holds. This completes the proof of Theorem \ref{main result - compatibility}.

\end{document}